\documentclass[aos,MSNbibl,nameyear,seceqn,dvips]{arximspdf}

\doi{10.1214/12-AOS1026} 
\volume{40}
\issue{3}
\pubyear{2012}
\firstpage{1906}
\lastpage{1933}

\makeatletter
\newcommand{\bz}{\mathbf{z}}
\newcommand{\bb}{\mathbf{b}}
\newcommand{\Deltam}{\bolds{\Delta}}
\newcommand{\tbx}{{\tilde{\bx}}}
\newcommand{\bx}{\mathbf{x}}
\newcommand{\bzero}{\mathbf{0}}
\newcommand{\bone}{\mathbf{1}}
\newcommand{\bA}{\mathbf{A}}
\newcommand{\bff}{\mathbf{f}}
\newcommand{\bF}{\mathbf{F}}
\newcommand{\bba}{\mathbf{a}}
\newcommand{\bI}{\mathbf{I}}
\newcommand{\bW}{\mathbf{W}}
\newcommand{\hW}{{\hat{W}}}
\newcommand{\hF}{{\hat{F}}}
\newcommand{\bZ}{\mathbf{Z}}
\newcommand{\bh}{\mathbf{h}}
\newcommand{\bX}{\mathbf{X}}
\newcommand{\alpham}{\bolds{\alpha}}
\newcommand{\betam}{\bolds{\beta}}
\newcommand{\deltam}{\bolds{\delta}}
\newcommand{\Sigmam}{\bolds{\Sigma}}
\newcommand{\etam}{\bolds{\eta}}
\newcommand{\R}{\mathbb{R}}

\newcommand{\cC}{\mathcal{C}}
\newcommand{\cF}{\mathcal{F}}
\newcommand{\cH}{\mathcal{H}}
\newcommand{\hbF}{{\hat{\bF}}}
\newcommand{\hbW}{{\hat{\bW}}}
\newcommand{\tdelm}{{\tilde{\deltam}}}
\newcommand{\cdelm}{{\check{\deltam}}}
\newcommand{\hdelm}{{\hat{\deltam}}}
\newcommand{\tbX}{{\tilde{\bX}}}

\newtheorem{lemma}{Lemma}
\newtheorem{theorem}{Theorem}
\newproclaim{remark}{Remark}
\newtheorem{corollary}{Corollary}
\makeatother

\begin{document}
\begin{frontmatter}

\title{Flexible generalized varying coefficient regression~models}
\runtitle{Flexible varying coefficient models}

\begin{aug}
\author[a]{\fnms{Young K.} \snm{Lee}\thanksref{t1}\ead[label=e1]{youngklee@kangwon.ac.kr}},
\author[b]{\fnms{Enno} \snm{Mammen}\ead[label=e2]{emammen@rumms.uni-mannheim.de}}
\and
\author[c]{\fnms{Byeong U.} \snm{Park}\corref{}\thanksref{t2}\ead[label=e3]{bupark@stats.snu.ac.kr}}
\thankstext{t1}{Supported by the Basic Science Research Program through
the NRF funded by the Korea government (MEST) (No. 2010-0021396).}
\thankstext{t2}{Supported by the NRF grant funded by the Korea
government (MEST) (No. 2010-0017437).}
\runauthor{Y. K. Lee, E. Mammen and B. U. Park}
\affiliation{Kangwon National University, Universit\"at Mannheim  and~Seoul~National~University}
\address[a]{Y. K. Lee\\
Department of Statistics\\
Kangwon National University\\
Chuncheon 200-701\\
Korea\\
\printead{e1}}

\address[b]{E. Mammen\\
Department of Economics\\
University of Mannheim\\
68131 Mannheim, L7, 3-5\\
Germany\\
\printead{e2}}

\address[c]{B. U. Park\\
Department of Statistics\\
Seoul National University\\
Seoul 151-747\\
Korea\\
\printead{e3}}
\end{aug}

\received{\smonth{1} \syear{2012}}
\revised{\smonth{6} \syear{2012}}

%
\begin{abstract}
This paper studies a very flexible model that can be used widely
to analyze the relation between a response and multiple
covariates. The model is nonparametric, yet renders easy
interpretation for the effects of the covariates. The model
accommodates both continuous and discrete random variables for the
response and covariates. It is quite flexible to cover the
generalized varying coefficient models and the generalized
additive models as special cases. Under a weak condition we give a
general theorem that the problem of estimating the multivariate
mean function is equivalent to that of estimating its univariate
component functions. We discuss implications of the theorem for
sieve and penalized least squares estimators, and then investigate
the outcomes in full details for a kernel-type estimator. The
kernel estimator is given as a solution of a system of nonlinear
integral equations. We provide an iterative algorithm to solve the
system of equations and discuss the theoretical properties of the
estimator and the algorithm. Finally, we give simulation results.
\end{abstract}

%


\end{frontmatter}

\section{Introduction}\label{sec1}
The varying coefficient regression model, proposed by \citet{hastie1993}, and studied by \citet{yang},
Roca-Pardinas and Sperlich (\citeyear{roca}) and \citet{lee2012}, is known to be a useful tool for analyzing the relation
between a response and a multivariate covariate. For a response
$Y$ and covariates $\bX$ and $\bZ$, they assumed that the mean
regression function, $m(\mathbf{x},\mathbf{z}) \equiv E(Y | \bX=\mathbf{x}, \bZ
=\mathbf{z})$,
takes the form $m(\mathbf{x},\mathbf{z}) = x_1 f_1(z_1)+ \cdots+ x_{d} f_{d}
(z_{d})$ for some unknown univariate functions $f_j$. The model is
simple in structure, gives easy interpretation, and yet is
flexible since the dependence of the response variable on the
covariates is modeled in a nonparametric way. The major hurdle in
the practical application of this model is that one needs to pair
up each ``$X$-covariate'' with only one ``$Z$-covariate.'' Typically
this is not the case in practice. In principle, each $X$-covariate
may interact with any number of $Z$-covariates to explain the
variation in $Y$. Moreover, it
is often difficult to differentiate $X$-types from $Z$-types in a
group of the covariates.

In this paper we are concerned with quite a more flexible setting
than the usual varying coefficient model.\vspace*{1pt} Suppose that we are
given a group of covariates, $X_1, \ldots, X_D$. Let $\bX=(X_1,
\ldots, X_D)^\top$ and $d$ be an integer such that $d \le D$. The
model we are interested in assumes that there is a link function,
say $g$, such that the mean function $m(\mathbf{x})\equiv E(Y|\bX=\mathbf{x})$
satisfies
%
\begin{equation}
\label{model} g\bigl(m(\mathbf{x})\bigr) = x_1 \biggl(\sum
_{k \in I_1}f_{1k}(x_{k}) \biggr)+ \cdots +
x_{d} \biggl(\sum_{k \in I_{d}}f_{dk}(x_k)
\biggr),
\end{equation}
where the index sets $I_j \subset\{1, 2, \ldots, D\}$ are known,
and each $I_j$ does not include $j$. The covariates that enter
into one of the coefficient functions $f_{jk}$, that is, $X_k$ for $k
\in\mathcal{C} \equiv\bigcup_{j=1}^d I_j$ are of continuous type. For
simplicity, we assume that $X_k$ with $k \in\cC$ are supported on
the interval $[0,1]$. We allow some of the covariates $X_j$, for
$1 \le j \le d$, to be discrete random variables. We also allow
that $\cC$ and $\{1, 2, \ldots, d\}$ may have common indices. Let
$\cC_0 = \cC\cap\{1, 2, \ldots, d\}$. The case $\cC_0 =
\varnothing$, that is, $\cC=\{d+1, \ldots, D\}$, corresponds to the
situation where one can distinguish between two groups of
covariates, $\{X_1, \ldots, X_d\}$ and $\{Z_1, \ldots, Z_p\}
\equiv\{X_{d+1}, \ldots, X_{d+p}\}$ with $D=d+p$. In this case,
the model reduces to
%
\begin{equation}
\label{modelxz} g\bigl(m(\mathbf{x},\mathbf{z})\bigr) = x_1 \biggl(
\sum_{k \in I_1}f_{1k}(z_{k}) \biggr)+
\cdots+ x_{d} \biggl(\sum_{k \in I_{d}}f_{dk}(z_k)
\biggr),
\end{equation}
where $I_j \subset\{1, 2, \ldots, p\}$ are index sets of $Z_k$.
The latter model arises, for example, when one takes all $X_j$,
for $1 \le j \le d$, to be discrete covariates. The above model reduces
further to the
nonparametric generalized additive model of \citet{yu} if we take $d=1$ and $X_1 \equiv1$.

The functions $f_{jk}$ in the representation (\ref{model}) are not
identifiable. To see this, consider the case where $d=D=3,
I_1=\{2,3\}, I_2=\{3\}, I_3=\{2\}$ so that $\cC=\cC_0=\{2,3\}$. In
this case, $x_1[f_{12}(x_2)+f_{13}(x_3)]+x_2 f_{23}(x_3)+x_3
f_{32}(x_2) =  x_1[g_{12}(x_2)+g_{13}(x_3)]+x_2 g_{23}(x_3)+x_3
g_{32}(x_2)$, if $g_{12}(x_2)=f_{12}(x_2)+c$,
$g_{13}(x_3)=f_{13}(x_2)-c$, $g_{23}(x_3)=f_{23}(x_3)+x_3$ and
$g_{32}(x_2)=f_{32}(x_2)-x_2$ for some constant $c$. To make all
$f_{jk}$ identifiable, we use the following constraints:
%
\begin{eqnarray}
\label{constraint1} \int f_{jk}(x_k) w_k(x_k)
\,dx_k &=& 0,\qquad  k \in\mathcal{C}, 1 \le j \le d,
\nonumber
\\[-8pt]
\\[-8pt]
\nonumber
\int x_k f_{jk}(x_k) w_k(x_k)
\,dx_k &=& 0, \qquad j,k \in\mathcal{C}_0
\end{eqnarray}
for nonnegative weight functions $w_k$, where $\mathcal{C}_0 = \{1,
2, \ldots, d\} \cap\mathcal{C}$. With these constraints we may
rewrite model (\ref{model}) as
%
\begin{eqnarray}
\label{modelnormed}g\bigl(m(\mathbf{x})\bigr) = \sum
_{j=1}^d \alpha_j x_j +
\mathop{\sum_{ j <
k }}_{j,k \in
\mathcal{C}_0 } \alpha_{jk} x_j
x_k + \sum_{j=1}^d
x_j \biggl(\sum_{k
\in I_j}
f_{jk}(x_k) \biggr).
\end{eqnarray}

We think that our approach broadens the field of applications of
varying coefficient models essentially. The link function allows
us to have a discrete response. Furthermore, our model frees us
from the restrictive settings of the usual varying coefficient
model that one should differentiate between two types of
covariates, $X$- and $Z$-type as in (\ref{modelxz}), and that each
covariate appears in only one ``nonlinear interaction term.'' In
the case where $X_1, \ldots, X_q$, for some $q \le D-2$, are
discrete and the remaining covariates are continuous random
variables, our approach allows fitting the full model, where
$I_j=\{q+1, \ldots, D\}$ for all $1 \le j \le q$ and $I_j=\{q+1,
\ldots, j-1, j+1, \ldots, D\}$ for $q+1 \le j \le D$:
%
\begin{eqnarray}
\label{model3}  g\bigl(m(\mathbf{x})\bigr)&=&x_1 \Biggl(
\sum_{k=q+1}^D f_{1k}(x_{k})
\Biggr)+ \cdots+x_{q} \Biggl(\sum_{k=q+1}^{D}
f_{pk}(x_{k}) \Biggr)
\nonumber
\\[-8pt]
\\[-8pt]
\nonumber
&&{}  + x_{q+1} \Biggl(\sum_{k=q+2}^D
f_{q+1,k}(x_{k}) \Biggr)+ \cdots+x_{D} \Biggl(\sum
_{k=q+1}^{D-1} f_{Dk}(x_{k})
\Biggr).
\end{eqnarray}
One may fit the full model in an exploratory analysis of the data,
and find a parsimonious model, deleting some of the functions
$f_{jk}$ in the full model, that fits the data well.

We stress that fitting model (\ref{model}) is not more complex than
fitting other varying coefficient models such as $g(m(\mathbf{x})) =
x_1f_{1}(x_{d+1})+ \cdots+ x_{d}f_{d}(x_{2d})$. The complexity is only
in notation and theory, yet it gives full flexibility in modeling via
varying coefficients.
Think of the case where the true model is $g(m(\mathbf{x}))=x_1 f_{12}(x_2)+
x_2f_{21}(x_1)$. In this case, $g(m(\mathbf{x}))$ may not be well approximated
by either $x_1f_{12}(x_2)$ or $x_2f_{21}(x_1)$ alone. Each additive
term in $g(m(\mathbf{x}))$ is interpreted as a $\operatorname{(linear)}\times
\operatorname
{(nonlinear)}$ interaction. With $x_2$ being held fixed, modeling by
$x_1 f_{12}(x_2)$ alone, for example, reflects only the linear effect
of $X_1$, while modeling by $x_1 f_{12}(x_2)+ x_2f_{21}(x_1)$
accommodates the nonlinear effect of $X_1$ as well.

\citet{xue} discussed a special case of model
(\ref{modelxz}) where one can differentiate between $X$-type and
$Z$-type variables and all $I_j=\{1,2, \ldots, p\}$. They proved
that the functions $f_{jk}$ with the constraints
(\ref{constraint1}) are identifiable. The essential assumption
was that the smallest eigenvalue of $E(\bX\bX^\top
|\bZ=\mathbf{z})$ is bounded away from zero where $\bX=(X_1, \ldots,
X_d)^\top$ and $\bZ=(Z_1, \ldots, Z_p)^\top$, although they put a
stronger one; see their condition (C2). Their approach cannot be
extended to our model~(\ref{model}). To see this, consider the
model~(\ref{model3}) and think of $E(\bX\bX^\top| \bX^c=\mathbf{x}^c)$
where $\bX= (X_1, \ldots, X_D)^\top$ and $\bX^c=(X_{q+1}, \ldots,\break
X_D)^\top$. Certainly, the matrix is singular for all $\mathbf{x}^c$ if
$D-q \ge2$. One of our main tasks in this paper is to relax the
requirement that the smallest eigenvalue of $E(\bX\bX^\top|
\bX^c=\mathbf{x}^c)$ is bounded away from zero; see assumption (A0)
in Section~\ref{sec2}. This weaker condition is typically satisfied. In
Lemma~\ref{idenlem} in the \hyperref[app]{Appendix}, we show that, under the weaker condition,
the $L_2$-norms of~$m$ and of the function tuple $(f_{jk}\dvtx  k \in
I_j, 1 \le j \le d)$ in the model~(\ref{model}) are equivalent,
modulo the norming constraints
(\ref{constraint1}). This has important
implications in estimating the model (\ref{model}). First, it
implies that the functions~$f_{jk}$ are identifiable. For some
types of estimators of $m$, it also gives directly the first-order
properties of the corresponding estimators of $f_{jk}$. In the
next section, we discuss the implications. In the subsequent
sections, we focus on the smooth backfitting method of estimating
the model (\ref{model}), where we also use the main idea contained
in Lemma~\ref{idenlem} to derive its asymptotic properties. The original idea
of smooth backfitting was introduced by \citet{mammen1999} for fitting additive models, and it has been further
developed in other contexts; see Yu, Park and Mammen (\citeyear{yu}) and Lee, Mammen and Park (\citeyear{lee2010,lee2012}), for example.

Earlier works on varying coefficient models focused on the model
$m(\mathbf{x},\mathbf{z})=x_1f_1(\mathbf{z})+ \cdots+x_{d}f_{d}(\mathbf{z})$,
where a single
covariate $\bZ$ (univariate or multivariate) enters into all
coefficient functions.
This model was proposed and studied by \citet{chen1993},
\citet{kauermann}, Fan and Zhang (\citeyear{fan1999}, \citeyear{fan2000}),
\citet{cai2000a}, \citet{cai2000b} and \citet{fan2003}. \citet{mammen2003} added a link function to this
model: $g(m(\mathbf{x},\mathbf{z}))=x_1f_1(\mathbf{z})+ \cdots+x_{d}f_{d}(\mathbf{z})$. The case
where $f_j$ are time-varying was treated by \citet{hoover},
Huang, Wu and Zhou (\citeyear{huang2002}, \citeyear{huang2004}), \citet{wang}, and
\citet{noh}. \citet{heim} also considered the case where all coefficient functions
are defined on a single 3D spatial space. Fitting these models is
simple. A univariate or multivariate smoothing across the single
variable $\bZ$, or on a time scale, or on a multidimensional spatial
space, gives directly estimators of $f_j$ without further
projection (by marginal integration or backfitting, e.g.)
onto relevant function spaces. However, this suffers from the curse of
dimensionality when the dimension of $\bZ$ is high. For this reason,
most works were focused on univariate $Z$.
Some time series models related to
the functional coefficient model, with $X_j$
being unobserved common factors that depend on time, have been
proposed and studied by, for example, Fengler, H\"ardle and Mammen (\citeyear{fengler})
and \citet{park}.

\section{Equivalence in entropies of function classes}\label{sec2}

In this section we will show that the nonparametric components,
$f_{jk}\dvtx  k \in I_j,  1 \le j \le d$, of our mod\-el~(\ref{modelnormed}) with constraints
(\ref{constraint1}) can be estimated with a
one-dimensional nonparametric rate. This means that our model
avoids the curse of dimensionality. It is easy to check that the
function $m(\mathbf{x})$ can be estimated with a one-dimensional rate.
This follows by application of results from empirical process
theory; see below. We will use that the $L_2$-norms of $m$ and of
the tuples $(\alpham; f_{jk}\dvtx  1 \leq j \leq d, k \in I_j)$ are
equivalent; see Lemma~\ref{idenlem} in the proofs section. Here,
$\alpham$ denotes the vector with elements $\alpha_{j}\dvtx 1 \leq j
\leq d$, $\alpha_{jk}\dvtx  j<k, j,k \in\mathcal{C}_0$. Our next result
uses this fact to show that the rate for the estimation of $m$
carries over to the estimation of $( \alpham; f_{jk}\dvtx  1 \leq j
\leq d, k \in I_j)$.

In the description of our method and in our theory we will also
make use of a different representation of the model
(\ref{modelnormed}). In this representation of the model, we
collect those coefficients that are functions of the same
continuous covariate and put them together as an additive
component. Suppose that, among $X_{1}, \ldots, X_d$ in model
(\ref{modelnormed}), there are $r$ ($0 \le r \le d$) variables
whose indices do not enter into $\cC$. Without loss of generality,
we denote them by $X_{1}, \ldots, X_{r}$. Let $p=D-r \ge2$ be the
number of indices in $\cC$. Thus $\cC=\{r+1, \ldots, r+p\}$ and
$\cC_0=\{r+1, \ldots, d\}$. Define
%
\begin{equation}
\label{xtild} \tbX_k =\{X_j\dvtx  r+k \in I_j,
1 \le j \le d\},\qquad 1 \le k \le p.
\end{equation}
The vector $\tbX_k$ is the collection of all $X_j$, for $1 \le j
\le d$, that have interactions with $X_{r+k}$ in the form of $X_j
f_{j,r+k}(X_{r+k})$. Thus, $\tbX_k$ does not include $X_{r+k}$ as
its element. Let $d_k$ denote the number of the index sets $I_j$
that contain $r+k$. Thus, $\tbX_k$ is of $d_k$-dimension.
Likewise, for a given vector $\mathbf{x}$, we denote the above
rearrangements of $\mathbf{x}$ by ${\tilde{\mathbf{x}}}_k,  1 \le k \le p$.
Also, define
$\bff_k =\{f_{j,r+k}\dvtx  r+k \in I_j,  1 \le j \le d\}$ for $1 \le k
\le p$. Then model (\ref{modelnormed}) can be represented as
%
\begin{equation}
\label{model1a} \qquad g\bigl(m(\mathbf{x})\bigr) =  \sum_{j=1}^d
\alpha_j x_j + \mathop{\sum_{ j < k }}_{
j,k \in
\mathcal{C}_0 }
\alpha_{jk} x_j x_k
+ {\tilde{\mathbf{x}}}_1^\top\bff_1(x_{r+1})
+ \cdots+ {\tilde{\mathbf{x}}}_p^\top\bff_p(x_{r+p}).
\end{equation}

To give an example of the above representation, consider the case
where $d=D=3$, $I_1=\{2,3\}$, $I_2=\{3\}$, $I_3=\{2\}$ so that
$\cC=\cC_0=\{2,3\}$. In this case, $r=1$, ${\tilde{\mathbf{x}}}_1 =(x_1,
x_3)^\top$, ${\tilde{\mathbf{x}}}_2=(x_1,x_2)^\top$, $\bff_1=(f_{12},
f_{32})^\top$, $\bff_2 = (f_{13}, f_{23})^\top$, and thus
\[
x_1\bigl[f_{12}(x_2)+f_{13}(x_3)
\bigr]+x_2 f_{23}(x_3)+x_3
f_{32}(x_2) = {\tilde{\mathbf{x}}}_1^\top
\bff_1(x_2) + {\tilde{\mathbf{x}}}_2^\top
\bff_2(x_3).
\]

Suppose now that one has an estimator $\hat m$ of $m$ with
%
\begin{equation}
\label{estmodmodel} g\bigl(\hat m(\mathbf{x})\bigr) = \sum
_{j=1}^d \hat\alpha_j x_j
+ \mathop{\sum_{ j < k}}_{ j,k \in\mathcal{C}_0 }\hat\alpha_{jk}
x_j x_k + \sum_{j=1}^d
x_j \biggl(\sum_{k \in I_j}\hat
f_{jk}(x_k) \biggr),
\end{equation}
where the estimated functions $\hat f_{jk}$ satisfy the
constraints in (\ref{constraint1}). We make the
following assumption:
\begin{longlist}
\item[(A0)] It holds that the product measure $\prod_{j=1}^D
P_{X_j}$ has a density w.r.t. the distribution $ P_{\bX}$ of
$\bX$ that is bounded away from zero and infinity on the support
of~$ P_{\bX}$. Here, $P_{X_j}$ is the marginal distribution of
$X_j$. The marginal distributions are absolutely continuous
w.r.t. Lebesgue measure or they are discrete measures with a
finite support. Furthermore, the weight functions $w_j$ for $j \in
\cC$ in the constraints in (\ref{constraint1}) are chosen so that $w_j/
p_{X_j}$ is bounded away from zero and infinity on the support of $
P_{X_j}$. Here, $p_{X_j}$ is the
density of $X_j$. The smallest eigenvalues of the matrices
$E[\tbX_k \tbX_k^\top| X_{r+k}=z_k]$ for $1 \le k \le p$ are
bounded away from zero for $z_k$ in the support of $p_{X_{r+k}}$.
\end{longlist}

The condition on $E[\tbX_k \tbX_k^\top| X_{r+k}=z_k]$ in (A0) is
typically satisfied. For example, consider the model
(\ref{model3}), where $X_1, \ldots, X_q$ are discrete random
variables whose indices do not enter into $\cC$. Thus, $r=q$ and
$p=D-q$. According to configuration (\ref{xtild}), we get
$\tbX_k =(X_j\dvtx  1 \le j \le D, j \neq q+k)^\top$
which is the covariate vector $\bX$ without $X_{q+k}$. In this
case, $E[\tbX_k \tbX_k^\top| X_{r+k}=z_k]$ is positive definite
if the support of the conditional distribution $P_{\tbX_k
|X_{q+k}=z_k}$ contains $(D-1)$ linearly independent vectors.

We get the following theorem for the rate of convergence of the
components of~$\hat m$.

\begin{theorem}\label{rateproposition} Suppose that assumption \textup{(A0)}
applies, and that an estimator $\hat m$ of $m$ with (\ref{estmodmodel})
and (\ref{constraint1}) satisfies that, for a null sequence $\kappa_n $,
%
\begin{equation}
\label{estrate} \int\bigl[g\bigl(\hat m(\mathbf{x})\bigr)-g\bigl(m(\mathbf{x})
\bigr)\bigr]^2 P_{\bX}(d\mathbf{x}) = O_p\bigl(
\kappa_n^2\bigr).
\end{equation}
Then it holds that
\[
\int\bigl[\hat f_{jk}(x_k) - f_{jk}(x_k)
\bigr]^2 p_{X_k}(x_k) \,d x_k=O_p
\bigl(\kappa_n^2\bigr)
\]
for all $k \in I_j,   1 \le j \le d$.
\end{theorem}

It is easy to construct estimators that fulfill (\ref{estrate}).
We will discuss this for the case where i.i.d. observations
$(\bX^i, Y^i)$ are made on the random vector $(\bX,Y)$ with $\bX^i
\equiv(X_1^i, \ldots, X_{D}^i)^\top$ of $D$-dimension. Examples
for estimators that fulfill (\ref{estrate}) are sieve estimators
or penalized least squares estimators. If one makes entropy
conditions on some function classes $\mathcal{F}_{jk}$, then Theorem
\ref{rateproposition} can be used to show that the entropy
conditions carry over to the class $\mathcal{M} =\{ m\dvtx  g(m(\mathbf{x}))
\mbox{ has the structure (\ref{modelnormed}) for some } \alpham
\in\mathcal{A} \mbox{ and } f_{jk} \in\mathcal{F}_{jk}, k \in
I_j,\break
1 \le j \le d\}$ for some compact set $\mathcal{A}$. Using empirical
process methods one can then show that sieve estimators or
penalized least squares estimators fulfill (\ref{estrate}). Below
we outline this for the case where $\mathcal{F}_{jk}$ are the classes
of $l$-times differentiable functions for some $l \ge2$.\vadjust{\goodbreak}

The penalized least squares estimator $(\hat{\alpham}^\mathrm{PEN};
\hat f_{jk}^\mathrm{PEN}\dvtx  1\le j \le d, k \in I_j)$ minimizes
%
\begin{eqnarray}
\label{eq2pen} &&n^{-1} \sum_{i=1}^n
\Biggl\{Y_i - g^{-1} \Biggl(\sum
_{j=1}^d \alpha_j X_j^i
+ \mathop{\sum_{ j < k }}_{ j,k \in\mathcal{C}_0
}\alpha_{jk}
X_j^i X_k^i
\nonumber
\\[-8pt]
\\[-8pt]
\nonumber
&&\hspace*{110pt}{} + \sum_{j=1}^d X_j^i
\sum_{k \in
I_j}f_{jk}\bigl(X_k^i
\bigr) \Biggr) \Biggr\}^{2} + \lambda_{n}^{2} J
(\bff),
\end{eqnarray}
where we put $J (\bff) = \sum_{k \in I_1} \int D_z^l f_{1k}(z)^2
\,dz + \cdots+ \sum_{k \in I_d} \int D_z^l f_{dk}(z)^2\,dz.$ Here,
the functions $\hat f_{jk}^\mathrm{PEN}$ are chosen so that
(\ref{constraint1}) holds, and for a function
$g$, $D_z^l g$ denotes its $l$th order derivative. We get the
following result for the rate of convergence of the penalized
least squares estimators $ \hat f_{jk}^\mathrm{PEN}$.

\begin{corollary}\label{corratepen} Suppose that all assumptions of
Theorem~\ref{rateproposition} hold, that the link function $g$ has an
absolutely bounded derivative, and that the estimators $\hat
f_{jk}^\mathrm{PEN}$ of $f_{jk}$ are defined by (\ref{eq2pen}). Suppose
that, for $k \in I_j, 1 \leq j \leq d$, the functions $ f_{jk}$ have
square integrable derivatives of order $l$. Furthermore, assume that
the (conditional) distribution of $\varepsilon_i = Y_i - m(\bX^i),   1
\le i \le n$, has subexponential tails. That is, there are constants
$t_{\varepsilon}, c_{\varepsilon} > 0$ such that
\[
\sup_{1 \leq i \leq n}E\bigl[ \exp\bigl( t |\varepsilon_{i}|\bigr) |
\bX^1,\ldots,\bX^n\bigr] < c_{\varepsilon}
\]
almost surely for $|t| \leq t_{\varepsilon}$. Choose $\lambda_{n}$
such that $\lambda_{n}^{-1} = O_{p}(n^{l/(2l+1)})$ and
$\lambda_{n} = O_{p}(n^{-l/(2l+1)})$. Then it holds that
\begin{eqnarray*}
\int\bigl[\hat f_{jk}^\mathrm{PEN}(z) - f_{jk}(z)
\bigr]^2 p_{X_k}(z) \,dz&=&O_p
\bigl(n^{-2l/(2l+1)}\bigr),
\\
\int D_z^l\hat f_{jk}^\mathrm{PEN}
(z)^2 \,dz&=&O_{p}(1)
\end{eqnarray*}
for all $k \in I_j,   1 \le j \le d$.
\end{corollary}

We now discuss sieve estimation of $m$. We will do this for spline
sieves. Denote by $\mathcal{G}_{n,c}$ the space of all spline
functions of order $l$ with knot points~$0$, $L_n^{-1}, 2 L_n^{-1},
\ldots, 1$ and with $l$th derivative absolutely bounded by $c$.
The spline sieve estimator $(\hat{\alpham}^\mathrm{SIEVE}; \hat
f_{jk}^\mathrm{SIEVE}\dvtx  1\le j \le d, k \in I_j)$ minimizes
%
\begin{eqnarray}
\label{eq2sieve}
n^{-1} \sum_{i=1}^n
\Biggl\{\!Y_i - g^{-1}\! \Biggl(\sum
_{j=1}^d \alpha_j X_j^i
+ \!\mathop{\sum_{ j < k }}_{ j,k \in\mathcal{C}_0
}\!\alpha_{jk}
X_j^i X_k^i+ \!\sum_{j=1}^d X_j^i
\sum_{k \in
I_j}f_{jk}\bigl(X_k^i
\bigr)\! \Biggr)\! \Biggr\}^{2}\hspace*{-40pt}
\end{eqnarray}
over all functions $ f_{jk}$ in $\mathcal{G}_{n,c}$. Again, the
functions $\hat f_{jk}^\mathrm{SIEVE}$ are chosen so that~(\ref{constraint1}) holds.\vspace*{2pt} We get the
following result for the rate of convergence of the sieve
estimators $ \hat f_{jk}^\mathrm{SIEVE}$.

\begin{corollary}\label{corratesieve} Suppose that all assumptions of
Theorem~\ref{rateproposition} hold, that the link function $g$ has an
absolutely bounded derivative and that the estimators $\hat f_{jk}^\mathrm{SIEVE}$ of $f_{jk}$ are defined by (\ref{eq2sieve}). Suppose that, for
$k \in I_j, 1 \leq j \leq d$, the functions $ f_{jk}$ have derivatives
of order $l$ that are absolutely bounded by $c$.
Furthermore, assume that $E | \varepsilon|^{\gamma} < \infty$
holds for some $\gamma> 2+ l^{-1}$. Choose $L_n$ such that
$L_n^{-1} = O(n^{-1/(2l+1)})$ and $L_n = O(n^{1/(2l+1)})$. Then
it holds that
\[
\int\bigl[\hat f_{jk}^\mathrm{SIEVE}(z) - f_{jk}(z)
\bigr]^2 p_{X_k}(z) \,dz=O_p
\bigl(n^{-2l/(2l+1)}\bigr)
\]
for all $k \in I_j,   1 \le j \le d$.
\end{corollary}

Both results, Corollaries~\ref{corratepen} and~\ref
{corratesieve}, can be generalized to quasi-likelihood estimation.
Then, the estimators are defined as in (\ref{eq2pen}) or
(\ref{eq2sieve}), respectively, but with the squared error
$(y-g^{-1}(u))^2$ replaced by $Q(g^{-1}(u),y)$. For the definition
of~$Q$, see the next section. One can show that the results of
Corollaries~\ref{corratepen} and~\ref{corratesieve} still hold
under the additional assumptions (A1) and (A2) of \citet{mammen1997}.
This can be proved along the lines of arguments of
the latter paper. In the subsequent sections, we discuss kernel
estimation of the model (\ref{model}).

\section{Estimation based on kernel smoothing}\label{theory}

We will introduce a kernel estimator based on backfitting
and develop a complete asymptotic theory for this
estimator. Again, we will do this for the case where i.i.d.
observations $(\bX^i, Y^i)$ are made on the random vector
$(\bX,Y)$. Model (\ref{model}) can be rewritten as (\ref
{modelnormed}) with
constraints (\ref{constraint1}). It can be
shown that the parameters $\alpha_j$ and $\alpha_{jk}$ in
(\ref{modelnormed}) can be estimated at a faster rate than the
nonparametric functions $f_{jk}$. For this reason we neglect the
parametric parts for simplicity of presentation. Thus we consider
model~(\ref{model}) with the constraints~(\ref{constraint1}). In this
setting, our alternative representation~(\ref{model1a}) becomes
%
\begin{equation}
\label{model1} g\bigl(m(\mathbf{x})\bigr) = {\tilde{\mathbf{x}}}_1^\top
\bff_1(x_{r+1}) + \cdots+ {\tilde{\mathbf{x}}}_p^\top
\bff_p(x_{r+p}).
\end{equation}

Let $Q$ be the quasi-likelihood function such that $\partial
Q(\mu,y)/\partial\mu= (y-\mu)/\break V(\mu)$, where $V$ is a function
for modeling the conditional variance $v(\mathbf{x})\equiv\operatorname{
var}(Y|\bX=\mathbf{x})$ by $v(\mathbf{x}) = V(m(\mathbf{x}))$. The
quasi-likelihood for
the mean regression function $m$ is then given by $\sum_{i=1}^n
Q(m(\bX^i), Y^i)$, and taking into account the structure of the
model (\ref{model1}) the quasi-likelihood for $\bff_k$ is
%
\begin{equation}
\label{fullQL} \sum_{i=1}^n Q
\bigl(g^{-1}\bigl(\tbX_1^{i\top} \bff_1
\bigl(X_{r+1}^{i}\bigr) + \cdots +\tbX_p^{i\top}
\bff_p\bigl(X_{r+p}^{i}\bigr)\bigr),
Y^i\bigr).
\end{equation}

We take the smooth backfitting approach [Mammen, Linton and Nielsen
(\citeyear{mammen1999}), Yu, Park and Mammen (\citeyear{yu}), Lee, Mammen and Park (\citeyear{lee2010,lee2012})]. We
maximize the integrated kernel-weighted quasi-likelihood
%
\begin{eqnarray}
\label{objective}  L_Q(\etam)& \equiv&\int
n^{-1} \sum_{i=1}^n Q
\bigl(g^{-1}\bigl(\etam_1(z_1)^\top
\tbX_1^{i} + \cdots +\etam_p(z_p)^\top
\tbX_p^{i}\bigr), Y^i \bigr)
\nonumber
\\[-8pt]
\\[-8pt]
\nonumber
&&\hphantom{\int
n^{-1} \sum_{i=1}^n }{}  \times K_\bh\bigl(\bX^{c,i},\mathbf{z}\bigr)
\,d\mathbf{z}
\end{eqnarray}
over the tuple of functions $(\etam_1, \ldots, \etam_p)$, each
$\etam_k$ being a vector of univariate functions that satisfy the
constraints of (\ref{constraint1}), where
$\bX^{c,i}=(X_{r+1}^{i}, \ldots, X_{r+p}^{i})^\top$ and $\mathbf{z}=(z_1,
\ldots, z_p)^\top$. Here and throughout the paper, we label the
elements of a tuple~$\etam$ in such a way that
%
\begin{eqnarray}
\label{repre} \etam(\mathbf{z})&=&\bigl(\etam_{1}(z_1)^\top,
\ldots, \etam_{p}(z_p)^\top
\bigr)^\top,
\nonumber
\\[-8pt]
\\[-8pt]
\nonumber
\etam_k& =& \{\eta_{j,r+k}\dvtx  r+k \in I_j, 1 \le j
\le d\}
\end{eqnarray}
with $r$ introduced at the beginning of Section~\ref{sec2}, and with this
representation the constraints of (\ref{constraint1}) on the elements
of $\etam$ are
%
\begin{eqnarray}
\label{constraint3} \int\eta_{jl}(u) w_l(u) \,du &=& 0,\qquad  r+1
\le l \le r+p, 1 \le j \le d;
\nonumber
\\[-8pt]
\\[-8pt]
\nonumber
\int u \eta_{jl}(u) w_l(u) \,du& =& 0,\qquad  r+1 \le j,l \le d.
\end{eqnarray}

We take a kernel such that
%
\begin{equation}
\label{kernel} \int K_{h_j} (u,v) \,dv =1 \qquad\mbox{for all values of } u.
\end{equation}
This kernel can be constructed from the standard kernel of the
form $K_{h_j}(u-v)\equiv h_j^{-1}K((u-v)/h_j)$, where $K$ is a
symmetric nonnegative function, in such a way that $K_{h_j}(u,v) =
K_{h_j}(u-v)/\int K_{h_j}(u-w)  \,dw$. It was used in Mammen, Linton and
Nielsen (\citeyear{mammen1999}), Yu, Park and Mammen (\citeyear{yu}) and Lee,
Mammen and Park (\citeyear{lee2010,lee2012}), and will be used in our technical arguments.

The smooth backfitting estimators of $\bff_k$ in our model
(\ref{model1}) are ${\hat{\bff}}_k$ which maximize $L_Q$ at
(\ref{objective}). In Section~\ref{imple}, we detail an iterative procedure
to get the estimators. Here, we provide their theoretical
properties. We begin with some notational definitions. We let
$p_{\bX^c}$ denote the marginal density of $\bX^c = (X_{r+1},
\ldots, X_{r+p})^\top$. Also, we let $p_j$ and $p_{jk}$ be the
marginal densities of $X_{r+j}$ and $(X_{r+j},X_{r+k})$,
respectively. Define $Q_r(u,y) = \partial^r
Q(g^{-1}(u),y)/\partial u^r$ and
$\bW_{j}(\mathbf{z};\etam)= (\bW_{j1}(\mathbf{z};\etam), \ldots,
\bW_{jp}(\mathbf{z};\etam) )$, where
\[
\bW_{jk}(\mathbf{z};\etam) = - E \bigl[Q_2 \bigl(
\tbX^{\top}\etam \bigl(\bX^c\bigr), m(\bX) \bigr)
\tbX_j\tbX_k^{\top} | \bX^c=
\mathbf{z} \bigr] p_{\bX^c}(\mathbf{z})
\]
for $1 \le j, k \le p$, and $\tbX=(\tbX_1^\top, \ldots,
\tbX_p^\top)^\top$. Let
\[
\bW(\mathbf{z};\etam)= \bigl(\bW_1(\mathbf{z};\etam)^\top,
\ldots, \bW_p(\mathbf{z};\etam)^\top\bigr)^\top=
\bW(\mathbf{z};\etam)^\top.
\]
Throughout this paper we write $\bW(\mathbf{z})=\bW(\mathbf{z};\bff)$, where
$\bff$ is the true tuple of the coefficient functions. With slight
abuse of notation, we also write
\begin{eqnarray*}
\bW_{jk}(z_j, z_k;\etam)& =& - E
\bigl[Q_2 \bigl(\tbX^{\top
}\etam\bigl(\bX^c\bigr),
m(\bX) \bigr)\tbX_j\tbX_k^{\top} |
X_{r+j}=z_j, X_{r+k}=z_k \bigr]
\\
&&{}  \times p_{jk}(z_j,z_k),
\\
\bW_{jj}(z_j;\etam) &=& - E \bigl[Q_2 \bigl(
\tbX^{\top}\etam\bigl(\bX^c\bigr), m(\bX) \bigr)
\tbX_j\tbX_k^{\top} | X_{r+j}=z_j
\bigr] p_j(z_j).
\end{eqnarray*}
It follows that
\begin{eqnarray*}
\bW_{jj}(z_j;\etam) &=& \int\bW_{jj}(
\mathbf{z};\etam)\, d \mathbf{z}_{-j},  \bW_{jk}(z_j,z_k;
\etam) = \int\bW_{jk}(\mathbf{z};\etam)\, d \mathbf{z}_{-(j,k)},
\\
\bW_{jj}(z_j) &= &\int\bW_{jj}(\mathbf{z})\, d
\mathbf{z}_{-j}, \bW_{jk}(z_j,z_k)
= \int\bW_{jk}(\mathbf{z})\, d \mathbf{z}_{-(j,k)}.
\end{eqnarray*}
Here and throughout the paper, $\mathbf{z}_{-j}$ for a given vector $\mathbf{z}$
denotes the vector without its $j$th entry, and $\mathbf{z}_{-(j,k)}$
without its $j$th and $k$th entries. Due to the conditions (A0)
and (A1), $\bW_{jj}(z_j)$ is positive definite for all $z_j \in
[0,1]$. However, $\bW(\mathbf{z})$ may not be positive definite.

The relevant space of functions we are dealing with is the one
that consists of tuples $\etam$ of univariate functions with the
representation at (\ref{repre}) such that its elements satisfy the
constraints (\ref{constraint3}). For any
two functions $\etam^{(1)}$ and~$\etam^{(2)}$ of this type, define
$\langle\etam^{(1)}, \etam^{(2)} \rangle_\# = \int
\etam^{(1)\top}(\mathbf{z}) \bW(\mathbf{z})\etam^{(2)}(\mathbf{z})  \,d\mathbf{z}$ whenever it
exists. We denote by $\cH(\bW)$ the resulting space of tuples
$\etam$. The space is equipped with the inner product $\langle
\cdot, \cdot\rangle_\#$. Let $\|\cdot\|_\#$ be its induced norm,
that is, $\|\etam\|_\#^2 = \int\etam(\mathbf{z})^\top\bW(\mathbf{z})
\etam(\mathbf{z})   \,d\mathbf{z}$.

Define $p_{j,\bX}^{(1)}(\mathbf{x})=\partial p_{\bX}(\mathbf{x})/\partial
x_{r+j}$. Likewise, define $m_{j}^{(1)}(\mathbf{x})=\partial
m(\mathbf{x})/\partial x_{r+j}$ and $m_{j}^{(2)}(\mathbf{x})=\partial^2
m(\mathbf{x})/\partial x_{r+j}^2$.
Let $\Deltam_{jk}$ denote a $d_j$-vector such that
its $\ell$th element $\Delta_{jk,\ell}=1$ if the $\ell$th element of $\tbx_j$ in
our rearrangement (\ref{xtild}) equals $x_{r+k}$ and $\Delta_{jk,\ell}=0$ otherwise. Define
a $(d_1+\cdots+d_p)$-vector $\Deltam_k$ by $\Deltam_k^\top=(\Deltam_{1k}^\top, \ldots, \Deltam_{pk}^\top)$.
In the assumptions given in the
\hyperref[app]{Appendix}, we assume $n^{1/5}h_j \rightarrow c_j$ as $n \rightarrow
\infty$ for some constants $0<c_j<\infty$. For such constants,
define
\[
\tilde {\betam}_j(z_j) = \bW_{jj}(z_j)^{-1} p_j(z_j) \sum_{k=1}^p c_k^2
E\bigl(\bb_{jk}(\bX) |X_{r+j}=z_j\bigr)\int t^2 K(t)\,dt,
\]
where $\bb_{jk}(\bX)$ are $d_j$-vectors given by
\begin{eqnarray*}
\bb_{jk}(\bX)&=&\biggl(m_k^{(1)}(\bX) - \frac{\Deltam_k^\top \bff(\bX^c)}{g'(m(\bX))}\biggr)\\
&&{}\times\biggl[\frac{\tbX_j}{V(m(\bX))g'(m(\bX))}
\frac{p_{k,\bX}^{(1)}(\bX)}{p_{\bX}(\bX)}\\
&&\quad\hphantom{\biggl[}{}-\tbX_j \Deltam_k^\top \bff(\bX^c)\biggl(\frac{V'(m(\bX))}{V(m(\bX))^2 g'(m(\bX))^2}+\frac{g''(m(\bX))}{V(m(\bX)) g'(m(\bX))^3} \biggr)\\
&&\hspace*{196pt}{}+ \frac{\Deltam_{jk}}{V(m(\bX))g'(m(\bX))}\biggr]\\
&&{}+ \frac{1}{2}\frac{\tbX_j}{V(m(\bX))g'(m(\bX))}\biggl(m_k^{(2)}(\bX) + \frac{g''(m(\bX))
(\Deltam_k^\top \bff(\bX^c))^2}{g'(m(\bX))^3}\biggr).
\end{eqnarray*}

Let $\betam_{*}(\mathbf{z})= (\betam_{*j}(z_j)\dvtx 1 \leq j \leq p )$ be a
solution of
%
\begin{eqnarray}
\label{biasterm} \betam_{*j}(z_j) = \tilde{
\betam}_j(z_j)- \sum_{k\neq j}
\int \bigl[\bW_{jj}(z_j)\bigr]^{-1}
\bW_{jk}(z_j,z_k) \betam_{*k}(z_k)
\,dz_k,
\nonumber
\\[-8pt]
\\[-8pt]
\eqntext{1 \le j \le p,}
\end{eqnarray}
and put $\betam_j(z_j)$ to be the normalized versions of
$\betam_{*j}(z_j)$ so that they satisfy the constraints
(\ref{constraint3}). Below in a theorem we
show that $\betam(\mathbf{z})\equiv(\betam_1(z_1)^\top, \ldots,\break
\betam_p(z_p)^\top )^\top$ is the asymptotic bias of the smooth
backfitting estimator ${\hat{\bff}}(\mathbf{z})= ({\hat{\bff
}}_1(z_1)^\top, \ldots,
{\hat{\bff}}_p(z_p)^\top )^\top$.

For the special case where the matrix $\bW(\mathbf{z})$ is invertible,
one has an interpretation of $\betam(\mathbf{z})$ as the projection of
the asymptotic bias of the full-dimensional local quasi-likelihood
estimator that maximizes
\[
n^{-1} \sum_{i=1}^n Q\bigl(g^{-1}\bigl(\etam_1(\bz)^\top \tbX_1^{i} +
\cdots +\etam_p(\bz)^\top\tbX_p^{i}\bigr),
Y^i\bigr)K_\bh(\bX^{c,i},\bz).
\]
One can check that its asymptotic bias, $\betam_{\mathrm{mlt}}(\bz)$,
is given by
\[
\betam_{ \mathrm{mlt}}(\bz) =\bW(\bz)^{-1}p_{\bX^c}(\bz)\sum_{k=1}^p c_k^2
E\bigl(\bb_k(\bX) | \bX^c=\bz\bigr) \int t^2 K(t)\, dt,
\]
where $\bb_k(\bX)^\top = (\bb_{1k}(\bX)^\top, \ldots, \bb_{pk}(\bX)^\top)$.
%
%

The tuple $\betam_\mathrm{mlt}$ does not belong to $\cH(\bW)$. The
asymptotic bias $\betam(\mathbf{z})$ of the smooth backfitting estimator
${\hat{\bff}}(\mathbf{z})$ is identical to the projection of~$\betam_\mathrm{mlt}$
onto $\cH(\bW)$, that is, $ \betam= \operatorname{argmin}_{\etam\in
\cH(\bW)}\|\betam_\mathrm{mlt}-\etam\|_\#$. This
projection interpretation of $\betam(\mathbf{z})$ is not available if
$\bW(\mathbf{z})$ is not invertible. In general, one has to define
$\betam(\mathbf{z})$ through the integral equations (\ref{biasterm}).

For the asymptotic variance of ${\hat{\bff}}$, we define
\begin{eqnarray*}
\Sigmam_j(z_j) &=&\frac{1}{c_j p_j(z_j)}\int
K^2(t) \,dt \biggl[E \biggl(\frac{\tbX_j \tbX_j^\top}{V(m(\bX))g'(m(\bX))^2} \Big| X_{r+j}=z_j
\biggr) \biggr]^{-1}
\\
&& {}\times E \biggl(\frac{v(\bX)\tbX_j\tbX_j^\top}{V(m(\bX
))^2g'(m(\bX
))^2} \Big| X_{r+j}=z_j
\biggr)
\\
&&{} \times \biggl[E \biggl(\frac{\tbX_j
\tbX_j^\top}{V(m(\bX))g'(m(\bX))^2} \Big| X_{r+j}=z_j
\biggr) \biggr]^{-1},
\end{eqnarray*}
where $v(\mathbf{x})$ is the conditional variance of $Y$ given $\bX=\mathbf{x}$.

\begin{theorem}\label{Thmrate}
Under \textup{(A0)} in Section~\ref{sec2} and those \textup{(A1)--(A5)} in the \hyperref[app]{Appendix},
there exists a unique maximizer ${\hat{\bff}}$ of the integrated
kernel-weighted quasi-likelihood (\ref{objective}) with
probability tending to one. The maximizer ${\hat{\bff}}$ satisfies
\begin{eqnarray*}
\int\bigl|{\hat{\bff}}(\mathbf{z}) - \bff(\mathbf{z})\bigr|^2
p_\bZ(\mathbf{z}) \,d\mathbf{z}&=& O_p
\bigl(n^{-4/5}\bigr),
\\
\sup_{z_j \in[2h_j, 1-2h_j]} \bigl|{\hat{\bff}}_j(z_j) -
\bff_j(z_j)\bigr| &=& O_p\bigl(n^{-2/5}
\sqrt{\log n}\bigr),
\end{eqnarray*}
where $|\cdot|$ denote the Euclidean norm.
\end{theorem}

\begin{theorem}\label{Thmdist}
Assume that \textup{(A0)} in Section~\ref{sec2} and those \textup{(A1)--(A5)} in the \hyperref[app]{Appendix}
hold. Then, for all $\mathbf{z}$ in the interior of the support of
$p_{\bX^c}$, it follows that $n^{2/5}({\hat{\bff}}_j(z_j)-\bff_j(z_j))$
are jointly asymptotically normal with mean $(\betam_1(z_1)^\top,
\ldots, \betam_p(z_p)^\top)^\top$ and variance $\operatorname{diag}(\Sigmam_j(z_j))$.
\end{theorem}

In the special case where $m(\mathbf{x})=x_1f_1(x_{p+1})+\cdots+x_p
f_d(x_{2p})$ for $\mathbf{x}=(x_1, \ldots, x_{2p})^\top$, thus the link
$g$ is the identity function and $Q(\mu,y)=-(y-\mu)^2/2$, the
asymptotic bias and variance of ${\hat{\bff}}_j(z_j)$ stated in
Theorem~\ref{Thmdist} coincide with those in Theorem~2 of Lee, Mammen and Park (\citeyear{lee2012}).
This can be seen by noting that $\bW(\mathbf{z})$ is
invertible in this case, and that $V=g'=1$, $\tbX= (X_1, \ldots,
X_p)^\top$, $\bX^c = (X_{p+1}, \ldots, X_{2p})^\top$ and
$\bW(\mathbf{z})^{-1} E(\tbX X_j | \bX^c=\mathbf{z})p_{\bX^c}(\mathbf{z})=\bone_j$
where $\bone_j$ is the $p$-dimensional unit vector with the $j$th
entry being equal to one.

Theorem~\ref{Thmdist} can be also viewed as an extension of
Theorem~2 of Yu, Park and Mammen (\citeyear{yu}). In the latter work, smooth
backfitting for the additive model $g(m(\mathbf{z})) = f_1(z_1) + \cdots
+ f_p(z_p)$ for a link $g$ was considered. As we mentioned
earlier, model (\ref{model1}) reduces to the above model by
taking $d_k\equiv1$, $\tbX_k =X_k \equiv1$ for $1 \le k \le p$
and $r=p$. In this case, $\bW(\mathbf{z})$ is not invertible so that the
projection interpretation of $\betam(\mathbf{z})$ is not valid. If one
replaces $m(\mathbf{x})$ in the formula of $\tilde{\betam}_j(z_j)$ by
$g^{-1}(f(\mathbf{x}^c))\equiv g^{-1}(f_1(x_{r+1})+\cdots+f_p(x_{r+p}))$,
$\bW_{jk}$ by the corresponding quantities for the latter model,
which are
\begin{eqnarray*}
W_{jk}(z_j,z_k) &=& \int V
\bigl(g^{-1}\bigl(f(\mathbf{z})\bigr)\bigr)^{-1}g'
\bigl(g^{-1}\bigl(f(\mathbf{z} )\bigr)\bigr)^{-2}p_{\bX^c}(
\mathbf{z}) \,d\mathbf{z}_{-(j,k)},
\\
W_{jj}(z_j) &=& \int V\bigl(g^{-1}\bigl(f(
\mathbf{z})\bigr)\bigr)^{-1}g'\bigl(g^{-1}
\bigl(f(\mathbf{z})\bigr)\bigr)^{-2}p_{\bX^c}(\mathbf{z}) \,d
\mathbf{z}_{-j},
\end{eqnarray*}
then one can verify that the solution of the system of integral
equations in~(\ref{biasterm}) concides with the asymptotic bias in
Yu, Park and Mammen (\citeyear{yu}). The asymptotic variance $\Sigmam_j(z_j)$
given above also reduces to the one in Yu, Park and Mammen~(\citeyear{yu}).

\begin{remark}\label{re1}
In the case where the link $g$ in model (\ref{model1}) is the identity
(or a linear) function and when the covariates $X_j$ are independent,
one may apply marginally a kernel smoothing method to estimating each
coefficient function. To see this, suppose that all $\tbX_j$ contains
$1$ as their first entry and any entry of $\tbX_j$ does not equal to
any of $X_{r+k}$, $k \neq j$. Then, $E(Y|\tbX_j, X_{r+j})=\tbX_j^\top
[\bff_j(X_{r+j})+\mathbf{c}_j]$ for some constant vector $\mathbf{c}_j$. This
means that $\bff_j$ minimizes $E[Y-\tbX_j^\top\etam_j(X_{r+1})]^2$
over $\etam_j$ subject to a normalization. Thus, the marginal smoothing
that minimizes $n^{-1}\sum_{i=1}^n (Y^i-\etam_j^\top\tbX_j^i)^2
K_{h_j}(X_{r+j}^i,z_j)$ for each $j$ and each point $z_j$ gives a
consistent estimator of $\bff_j(z_j)$. This marginal smoothing approach
is not valid, even with independent covariates, in case the link
function $g$ is nonlinear. In the latter case, one needs to use a
projection method such as the smooth backfitting defined above, or a
marginal integration technique, to obtain appropriate estimators.
\end{remark}

\section{Implementation}\label{imple}

In this section, we discuss how to find ${\hat{\bff}}$ maximizing~$L_Q$ at
(\ref{objective}). Our method of finding ${\hat{\bff}}$ is based on an
iteration scheme. By considering the Fr\'echet differentials of
$L_Q$, we see that
%
\begin{eqnarray}
\label{esteq} \int n^{-1}\sum_{i=1}^n
Q_1 \bigl(\tbX^{i\top}{\hat{\bff}}(\mathbf{z}),
Y^i \bigr)\tbX_j^i K_\bh\bigl(
\bX^{c,i},\mathbf{z}\bigr) \,d\mathbf{z}_{-j} =
\bzero_j,
\nonumber
\\[-8pt]
\\[-8pt]
\eqntext{1 \le j \le p, z_j \in[0,1],}
\end{eqnarray}
where $\bzero_j$ is the zero vector of dimension $d_j$. The system
of equations is nonlinear. We take the Newton--Raphson\vspace*{1pt} approach to
find a solution by iteration. For a vector of functions
$(\etam_1^\top, \ldots, \etam_p^\top)^\top$ where $\etam_j(\mathbf{z}) =
\etam_j(z_j)$, define
%
\begin{eqnarray}
\label{def-F-j} \hbF_j(\etam) (z_j) = \int
n^{-1}\sum_{i=1}^n
Q_1 \bigl(\tbX^{i\top}\etam(\mathbf{z}), Y^i
\bigr)\tbX_j^i K_\bh\bigl(\bX^{c,i},
\mathbf{z}\bigr) \,d\mathbf{z}_{-j}.
\end{eqnarray}
The system of equations in (\ref{esteq}) is then expressed as
$\hbF_j({\hat{\bff}})=\bzero_j,   1 \le j \le p$. Our algorithm
runs an
\textit{outer} iteration which is based on a Newton--Raphson
approximation of the system of equations. Each outer-step solves a
linearized system of equations to update the approximation of
${\hat{\bff}}$, which requires an additional iteration, called \textit{inner}
iteration.

To describe the algorithm, suppose that we are at the $s$th
outer-step to update ${\hat{\bff}}^{[s-1]}$ in the previous outer-step.
Considering the Fr\'echet differentials\vadjust{\goodbreak} of $\hat\bF_j $ at
${\hat{\bff}}^{[s-1]}$, we get the following approximation: for $1
\le j
\le p$,
%
\begin{eqnarray}
\label{NR} \hat\bF_j ({\hat{\bff}}) (z_j)
&\simeq&\hat\bF_j \bigl({\hat{\bff}}^{[s-1]}\bigr)
(z_j)\nonumber\\[-2pt]
&&{}+ \sum_{k=1}^p \int
n^{-1}\sum_{i=1}^n
Q_2 \bigl(\tbX^{i\top}{\hat{\bff}}^{[s-1]}(\mathbf{z}),
Y^i \bigr)
\\[-2pt]
&&\hspace*{38pt}{} \times \tbX_j^i \tbX_k^{i\top}
\bigl[{\hat{\bff}}_{k}(z_k)-{\hat{\bff}}_{k}^{[s-1]}(z_k)
\bigr]K_\bh\bigl(\bX^{c,i},\mathbf{z}\bigr) \,d
\mathbf{z}_{-j}.\nonumber
\end{eqnarray}
This gives an updating equation for ${\hat{\bff}}^{[s]}$. Define, for $1
\le j \le p$ and for $1 \le j \neq k \le p$,
\begin{eqnarray*}
\hbW_{jk}^{[s]}(z_j,z_k) &=& - \int
n^{-1}\sum_{i=1}^n
Q_2 \bigl(\tbX^{i\top}{\hat{\bff}}^{[s]}(\mathbf{z}),
Y^i \bigr)\tbX_j^i \tbX_k^{i\top}
K_\bh\bigl(\bX^{c,i},\mathbf{z}\bigr) \,d\mathbf{z}_{-(j,k)},
\\[-2pt]
\hbW_{jj}^{[s]}(z_j) &=& - \int n^{-1}
\sum_{i=1}^n Q_2 \bigl(
\tbX^{i\top}{\hat{\bff}}^{[s]}(\mathbf{z}), Y^i \bigr)
\tbX_j^i \tbX_j^{i\top}
K_\bh\bigl(\bX^{c,i},\mathbf{z}\bigr) \,d\mathbf{z}_{-j}.
\end{eqnarray*}
Define $\tdelm_j^{[s]}(z_j) = [\hbW_{jj}^{[s]}(z_j)]^{-1}
\hat\bF_j ({\hat{\bff}}^{[s]})(z_j)$. Also, let $\hdelm_j^{[s]}(z_j) =
{\hat{\bff}}_j^{[s]}(z_j)-\break {\hat{\bff}}_j^{[s-1]}(z_j)$. We get from
(\ref{NR}) the
following linearized system of updating equations:
%
\begin{eqnarray}
\label{backeq} \qquad\hdelm_j^{[s]}(z_j) =
\tdelm_j^{[s-1]}(z_j) - \sum
_{k\neq j} \int \bigl[\hbW_{jj}^{[s-1]}(z_j)
\bigr]^{-1}\hbW_{jk}^{[s-1]}(z_j,z_k) \hdelm_k^{[s]}(z_k)
\,dz_k,
\nonumber
\\[-10pt]
\\[-10pt]
\eqntext{ 1 \le j \le p.}
\end{eqnarray}
Solving the system\vspace*{-2pt} of equations (\ref{backeq}) for
$\hdelm_j^{[s]}$ and then updating ${\hat{\bff}}_j^{[s-1]}$ by
${\hat{\bff}}_j^{[s]}={\hat{\bff}}_j^{[s-1]}+\hdelm_j^{[s]}$
constitutes the $s$th
step in the outer iteration. Below in Theorem~\ref{Thmalgo}, we
show that the outer iteration converges to a solution that
satisfies (\ref{esteq}).

The system\vspace*{-2pt} of equations (\ref{backeq}) cannot be solved since to
get $\hdelm_j^{[s]}$ one requires knowledge of the other $\hdelm_k^{[s]},
k \neq j$. To solve (\ref{backeq}) we need an \textit{inner}
iteration. Suppose that we are at the $\ell$th inner-step of the
$s$th outer-step to update $\hdelm_j^{[s,\ell-1]},   1 \le j \le
p$, in the previous inner iteration step. We apply~(\ref{backeq}):
%
\begin{eqnarray}
\label{innerit} \hdelm_j^{[s,\ell]}(z_j)& =&
\tdelm_j^{[s-1]}(z_j)- \sum
_{k\neq j} \int\bigl[\hbW_{jj}^{[s-1]}(z_j)
\bigr]^{-1}\hbW_{jk}^{[s-1]}(z_j,z_k)
\nonumber
\\[-11pt]
\\[-11pt]
\nonumber
&&\hspace*{82pt}{}\times
\hdelm_k^{[s,\ell-1]}(z_k)
\,dz_k,\qquad
 1 \le j \le p.
\end{eqnarray}
Existence of a unique solution of (\ref{backeq}) and the
convergence of the inner iteration to the solution are
demonstrated below in Theorem~\ref{Thmalgo}. For the starting
values $\hdelm_{j}^{[s,0]}$ in the inner iteration of the $s$th
outer-step, one may use the limit of the inner iteration in the
previous outer-step $\hdelm_{j}^{[s-1,\infty]}$.\vadjust{\goodbreak}

For a convergence criterion of the outer iteration, one may check
whether the values of the left-hand side of (\ref{esteq}) are
sufficiently small, or use the difference between the two updates
${\hat{\bff}}_k^{[s-1]}$ and ${\hat{\bff}}_k^{[s]}$:
%
\begin{equation}
\label{crit} \int\bigl|{\hat{\bff}}_{k}^{[s]}(z_k)- {
\hat{\bff}}_{k}^{[s-1]}(z_k)\bigr|^2
\,dz_k.
\end{equation}
In the latter case, one should use the normalized versions of the
updates. Recall the configuration of $\etam_k$ in (\ref{repre}).
The normalized version of a given set of tuples $\etam_{*k}$ may
be obtained by the following formulas. Let the weight functions
$w_l$ be normalized so that $\int w_l(u)  \,du=1$. Then
%
\begin{eqnarray}
\label{normali1} \eta_{jl}(u)&=& \eta_{jl*}(u)\nonumber\\
&&{}-\int
\eta_{*jl}(u)w_l(u) \,du,\qquad
1 \le j \le r \mbox{ or } d+1 \le l \le r+p,\\
\eta_{jl}(u)&=& \eta_{jl*}(u)-a_{jl}-b_{jl}
u,\qquad r+1 \le j, l \le d,\nonumber
\end{eqnarray}
where
\begin{eqnarray*}
a_{jl}&=& \int\eta_{*jl}(u)w_l(u) \,du -
b_{jl} \int u w_l(u) \,du,
\\
b_{jl}&=& \biggl[\int \biggl(u-\int t w_l(t) \,dt
\biggr)^2w_l(u) \,du \biggr]^{-1}
\\
&&{} \times \int \biggl(u-\int t w_l(t) \,dt \biggr)\eta_{*jl}(u)w_l(u)
\,du.
\end{eqnarray*}
One should also use the normalized $\hdelm_{jl}$ for the
convergence of the inner iteration.

\begin{theorem}\label{Thmalgo}
Assume that \textup{(A0)} in Section~\ref{sec2} and \textup{(A1)--(A5)} in the \hyperref[app]{Appendix}
hold. Then there exist constants $0<C_1,  \tau<\infty$ and
$0<\gamma<1$ such that, if the initial choice ${\hat{\bff}}^{[0]}$
satisfies
%
\begin{equation}
\label{inicond} \int\bigl|{\hat{\bff}}^{[0]}(\mathbf{z}) - {\hat{\bff}}(
\mathbf{z})\bigr|^2 p_{\bX^c}(\mathbf{z}) \,d\mathbf{z} \le
\tau^2
\end{equation}
with probability tending to one, then $\int|{\hat{\bff}}^{[s]}(\mathbf{z}) -
{\hat{\bff}}(\mathbf{z})|^2  p_{\bX^c}(\mathbf{z})  \,d\mathbf{z} \le  C_1
4^{-(s-1)}$
$\gamma^{2^s-1}$ with probability tending to one. Also, for each
outer-step there exists a solution of the system of equations
(\ref{backeq}) that is unique, and the inner iteration converges
at a geometric rate. Furthermore, if the initial choice
${\hat{\bff}}^{[0]}$ satisfies (\ref{inicond}) with probability
tending to
one, then there exist some constants $0<C_2<\infty$ and $0<\rho<1$
such that, with probability tending to one, $\int
|\hdelm^{[s,\ell]}(\mathbf{z}) - \hdelm^{[s,\infty]}(\mathbf{z})|^2
p_{\bX^c}(\mathbf{z})  \,d\mathbf{z} \le  C_2 \rho^{2\ell}$ for sufficiently
large $s$, where $\hdelm^{[s,\infty]}$ is a solution of the system
of equations in (\ref{backeq}).
\end{theorem}

The theorem shows that the number of iterations that is needed for
a~desired accuracy of the calculation of the backfitting estimator
does not depend on the sample size. If the desired accuracy is of
order $n^{-c}$ for some constant~$c$, then a logarithmic number of
iterations suffices. Thus the complexity of the algorithm only
increases very moderately for increasing sample size. We have no
good bound on the required accuracy of the starting values, that is,
on the choice of $\tau$. In our practical experience the algorithm
was very robust against poor choices of the starting value. In
fact, in the simulation study we chose ${\hat{\bff}}^{[0]}=\bzero$
and it
worked quite well. A more deliberate choice is a~version of the
marginal integration estimator studied by \citet{yang}, or a spline estimator that we discussed in Section~\ref{sec2}.
These are consistent so that they satisfy the condition
(\ref{inicond}), but they cost additional numerical computation.

As for the choices of the bandwidths $h_j$, one
may estimate the optimal bandwidths $h_j^\mathrm{opt}=c_j^* n^{-1/5}$,
where $\mathbf{c}^*=(c_1^*, \ldots, c_p^*)$ is defined by
%
\begin{equation}
\label{optband} \mathbf{c}^* = \arg\min_\mathbf{c} \sum
_{j=1}^p \int \bigl[ \bigl|\betam_j(z_j,
\mathbf{c})\bigr|^2 + \operatorname{ trace}\bigl(\Sigmam_j(z_j,c_j)
\bigr) \bigr]p_{Z_j}(z_j) \,dz_j.
\end{equation}
Here, we write $\betam_j(z_j,\mathbf{c})$ and $\Sigmam_j(z_j,c_j)$,
instead of $\betam_j(z_j)$ and $\Sigmam_j(z_j)$ as defined
in Section~\ref{theory}, to stress their dependence on the vector of the bandwidth
constants $\mathbf{c}=(c_1, \ldots, c_p)$. To describe a simple plug-in method,
get parametric estimates of $f_{jk}$ by maximizing (\ref{fullQL}) over
the class of $p$th order polynomials
$f_{jk}(x)=a_{jk}^{(0)}+a_{jk}^{(1)}x + \cdots+a_{jk}^{(p)}x^p$, and
obtain a kernel estimate of the density $p_\bX$. Then one can estimate
$\betam_j$ by plugging these estimates into the formulas of $\tilde
{\betam}_j,   \bW_{jj}, \bW_{jk}\ (j\neq k)$ and solving the system
of equations in (\ref{biasterm}) by iteration. One can also estimate
$\Sigmam_j$. Put these estimates of~$\betam_j$ and $\Sigmam_j$ into the
right-hand side of (\ref{optband}) to get an estimate of $\mathbf{c}^*$.
A~similar idea was adopted by Lee, Mammen and Park (\citeyear{lee2012}) for the classical varying
coefficient model. An alternative way is to develop a method similar to
the penalized least squares bandwidth selector proposed by \citet{mammen2005}. This would need higher-order stochastic expansions for the
quasi-likelihood of the smooth backfitting estimators. Finally, we want
to mention that $\betam_j$ depends on the whole vector $\bf c$ contrary
to $\Sigmam_j$, the latter only involving $c_j$. This is not the case
with the local linear smooth backfitting where both depend on $c_j$
only, see the next section. Thus, a grid search for~$\mathbf{c}^*$ at~(\ref
{optband}) may be computationally expensive for large $p$. In this
case, one may apply an iteration scheme which, in each iteration step,
updates~$c_j$ for $1 \le j \le p$ one by one with the other $c_k,   k
\neq j$, being held fixed at the values obtained in the previous step.

\section{Extension to higher-order local smoothing}\label{lpsbf}

In the previous two sections we considered smooth backfitting
based on Nadaraya--Watson smoothing. Here,\vadjust{\goodbreak} we discuss its extension
to local polynomial smoothing. We focus on the local linear case.
The extension to the general case is immediate, but needs more
involved notation. For a function $\eta$ of interest, the basic
idea of local linear smoothing is to approximate $\eta(u)$ for $u$
near a point $z$ by $\eta(z)+\eta'(z)(u-x)$, where $\eta'$ is the
first derivative of $\eta$. Thus, we maximize
\begin{eqnarray*}
&&\int n^{-1} \sum_{i=1}^n Q
\bigl(g^{-1}\bigl(\tbX_1^{i\top} \etam_1
\bigl(z_1,X_{r+1}^{i}\bigr) + \cdots +
\tbX_p^{i\top}\etam_p\bigl(z_p,X_{r+p}^{i}
\bigr)\bigr), Y^i \bigr)
\\[-2pt]
&&\hspace*{22pt}\qquad{} \times K_\bh\bigl(\bX^{c,i},\mathbf{z}\bigr)
\,d\mathbf{z},
\end{eqnarray*}
where $\etam_j(z_j,X_{r+j}^{i}) = \etam_{j0}(z_j) +
\etam_{j1}(z_j)h_j^{-1}(X_{r+j}^{i}-z_j)$. The maximizers, denoted
by ${\hat{\bff}}_{j0}$ and ${\hat{\bff}}_{j1}$ which correspond to
$\etam_{j0}$
and $\etam_{j1}$, respectively, are the estimators of the true
$\bff_j$ and $h_j\bff_j'$, where $\bff_j'$ is the vector of the
first derivatives of the entries in $\bff_j$. Again, ${\hat{\bff}}_{j0}$
should be normalized according to (\ref{normali1}).

To describe the algorithms, write ${\hat{\bff}}_{j}=({\hat{\bff
}}_{j0}^\top,
{\hat{\bff}}_{j1}^\top)^\top$. They satisfy $\hbF_j({\hat{\bff
}}_{1}, \ldots,\break
{\hat{\bff}}_{p})=\bzero_j$, $1 \le j \le p$, where $\bzero_j$
denotes now
the zero vector of dimension $2d_j$,
\begin{eqnarray*}
\hbF_j(\etam_{1}, \ldots, \etam_{p})& =&
\int n^{-1}\sum_{i=1}^n
Q_1 \Biggl(\sum_{j=1}^p
\tbX_j^{i\top}\etam_j\bigl(z_j,X_{r+j}^{i}
\bigr) , Y^i \Biggr)
\\[-2pt]
&&\hspace*{44pt}{} \times \bba\bigl(X_{r+j}^{i},z_j\bigr)
\otimes\tbX_j^i K_\bh\bigl(
\bX^{c,i},\mathbf{z}\bigr) \,d\mathbf{z}_{-j}
\end{eqnarray*}
$\etam_j =(\etam_{j0}^\top, \etam_{j1}^\top)^\top$ and
$\bba(X_{r+j}^{i},z_j)=(1,h_j^{-1}(X_{r+j}^{i}-z_j))^\top$. The
expressions for the updating equations at (\ref{backeq}) and
(\ref{innerit}) are unchanged if, writing
$\bA_{jk}^i=\bba(X_{r+j}^{i},z_j)\bba(X_{r+j}^{i},z_k)^\top$, we
redefine $\hbW_{jk}^{[s]}$, $1 \le j,k \le p$, by
\begin{eqnarray*}
\hbW_{jk}^{[s]}(z_j,z_k)& =& - \int
n^{-1}\sum_{i=1}^n
Q_2 \Biggl(\sum_{j=1}^p
\tbX_j^{i\top}{\hat{\bff}}_{j}^{[s]}
\bigl(z_j,X_{r+j}^{i}\bigr), Y^i
\Biggr)
\\[-2pt]
&&\hspace*{55pt}{} \times \bA_{jk}^i\otimes\tbX_j^i
\tbX_k^{i\top} K_\bh \bigl(\bX^{c,i},
\mathbf{z}\bigr) \,d\mathbf{z}_{-(j,k)},
\\[-2pt]
\hbW_{jj}^{[s]}(z_j) &=& - \int n^{-1}
\sum_{i=1}^n Q_2 \Biggl(\sum
_{j=1}^p \tbX_j^{i\top}{
\hat{\bff}}_{j}^{[s]}\bigl(z_j,X_{r+j}^{i}
\bigr), Y^i \Biggr)
\\[-2pt]
&&\hspace*{55pt}{} \times \bA_{jj}^i\otimes\tbX_j^i
\tbX_j^{i\top} K_\bh\bigl(\bX^{c,i},
\mathbf{z}\bigr) \,d\mathbf{z}_{-j}.
\end{eqnarray*}

Let $\Sigmam_j$ be defined as in Section~\ref{theory}, and define
$\betam_j^{LL}(z_j)$ to be the normalized versions of $c_j^2 \int t^2
K(t)   \,dt  \bff_j''(z_j)/2$ obtained by
(\ref{normali1}), where $\bff_j''$ is the vector
of the second derivatives of the entries in $\bff_j$.

\begin{theorem}
Under \textup{(A0)} in Section~\ref{sec2} and \textup{(A1)--(A5)} in the \hyperref[app]{Appendix},
Theorems~\ref{Thmrate} and~\ref{Thmalgo} remain valid for the
local linear smooth backfitting estimators $({\hat{\bff}}_{j0},
{\hat{\bff}}_{j1})$ and for their algorithms, respectively. As a
version of
Theorem~\ref{Thmdist}, $n^{2/5}({\hat{\bff}}_{j0}(z_j)-\bff_j(z_j))$ are
jointly asymptotically normal with mean $(\betam_1^{LL}(z_1)^\top,
\ldots, \betam_p^{LL}(z_p)^\top)^\top$ and variance $\operatorname{
diag}(\Sigmam_j(z_j))$.\vadjust{\goodbreak}
\end{theorem}

\section{Simulation study}\label{simst}

In the simulation study, we considered a binary response $Y$
taking values $0$ and $1$, and took the following model for the
mean function $m(\mathbf{x})$:
%
\begin{eqnarray}
\label{simodel} g\bigl(m(\bX)\bigr)&=& f_{02}(X_2) +
f_{03}(X_3) + X_1\bigl(f_{12}(X_2)+f_{13}(X_3)
\bigr)
\nonumber
\\[-8pt]
\\[-8pt]
\nonumber
&&{} + X_3 f_{32}(X_2)+X_2
f_{23}(X_3),
\end{eqnarray}
where $g(u)=\log(u/(1-u))$ is the logit link and $f_{02}(z)=z^2,
f_{03}(z)=4(z-0.5)^2,   f_{12}(z)=z,  f_{13}(z)=\cos(2\pi
z),  f_{32}(z)=e^{2z-1},  f_{23}(z)=\sin(2\pi z)$. The covariate
$X_1$ was a discrete random variable having $\operatorname{Bernoulli}(0.5)$
distribution, and $X_2$ and $X_3$ were $\operatorname{uniform} (0,1)$ random
variables. The three covariates were independent. We chose two
sample sizes $n=500$ and $1000$. The number of samples was $500$.
For the initial estimate, we used
${\hat{\bff}}^{[0]}=\bzero$. The weight functions $w_l$ were
$w_l(z)=I_{[0,1]}(z)$ for all $l$. We used the Epanechnikov kernel function
$K(u)=(3/4)(1-u^2)I_{[-1,1]}(u)$ and took the theoretically
optimal bandwidths as defined at (\ref{optband}), which were
$h_1^\mathrm{opt}=0.4328,   h_2^\mathrm{opt}=0.2789$ for $n=500$, and
$h_1^\mathrm{opt}=0.3768,   h_2^\mathrm{opt}=0.2428$ for $n=1000$, in
our simulation setting.

In the simulation, we also computed the cubic spline estimates
with $K$ knots placed evenly on the interval $[0,1]$. We used the
power basis for cubic splines: $s_0(z)=1, s_1(z)=z, s_2(z)=z^2,
s_3(z)=z^3, s_{3+k}(z)=(z-\xi_k)_+^3$, where $\xi_k$ are the knot
points. If one applies directly the power basis to the model
(\ref{simodel}), one may suffer from ``near singularity'' of the
resulting design matrix. This is because the functions $f_{jk}$
without satisfying our constraints are not identifiable. Taking
into consideration the constraints, we adjusted the power basis so
that $s_1$ is orthogonal to $s_0$, and $s_j$ for $2 \le j \le
K+3$, are orthogonal to $s_0$ and $s_1$. The dimension of the
power basis for the cubic spline approximation of the model
(\ref{simodel}) equals $6K+19$. The number of knots taken was
$K=1$ which gave the best performance. The performance of the
spline estimators got worse quickly as $K$ increased.

Table~\ref{tab1} shows the results based on 500 datasets. For each
component function $f_{jk}$, the table provides the integrated
mean squared error (IMSE), $\int E[\hat{f}_{jk}(z)-f_{jk}(z)]^2
\,dz$. The main lesson is that the spline estimators have much
larger variances than the smooth backfitting estimators, while the
former have smaller biases. Overall, the smooth backfitting method
works quite well. Comparing the values of IMSEs for the two sample
sizes, the results for the smooth backfitting method reflect the
asymptotic effects fairly well. Note that the theoretical
reduction of IMSE from $n=500$ to $n=1000$ equals $(0.5)^{4/5}
\simeq0.574$. In the simulation we also found the iterative
algorithm of the smooth backfitting method in Section~\ref{imple} converged
very fast. The outer loop typically converged in five iterations
with the criterion value $10^{-4}$ for the normalized difference
(\ref{crit}), and that the inner loop converged in three
iterations.

\begin{table}
\caption{Integrated mean squared errors (IMSE),
integrated squared biases (ISB) and integrated~variance~(IV)~of
the two methods, cubic spline (SPL) and
smooth~backfitting~(SBF),~for~the~model~(\protect\ref{simodel})}\label{tab1}
\begin{tabular*}{\textwidth}{@{\extracolsep{\fill}}lcccccccc@{}}
\hline
& & &$\bolds{f_{02}}$ &$\bolds{f_{12}}$&$\bolds{f_{32}}$ &$\bolds{f_{03}}$&$\bolds{f_{13}}$ &\multicolumn{1}{c@{}}{$\bolds{f_{23}}$}\\
\hline
SPL&$n=500$ &IMSE &0.2607 &0.2145 & 0.5034 &0.2631 &0.2274 & 0.5857\\
& & ISB &0.0013 &0.0004 &0.0011 &0.0006 &0.0044 &0.0159 \\
& & IV &0.2594 &0.2141 &0.5022 &0.2624 &0.2230 &0.5699 \\[3pt]
&$n=1000$ &IMSE & 0.1106 & 0.0817 & 0.2122 &0.1074 &0.0938 & 0.2453\\
& & ISB &0.0001 &0.0004 &0.0001 &0.0001 &0.0006 &0.0104 \\
& & IV &0.1105 &0.0813 &0.2121 &0.1073 &0.0932 &0.2349 \\[6pt]
SBF&$n=500$ &IMSE & 0.0315 & 0.0399 & 0.0274 &0.1071 &0.1073 & 0.1685\\
& & ISB &0.0035 &0.0112 &0.0128 &0.0090 &0.0543 &0.0808 \\
& & IV &0.0280 &0.0288 &0.0147 &0.0981 &0.0531 &0.0877 \\[3pt]
&$n=1000$ &IMSE & 0.0214 & 0.0210 & 0.0254 &0.0526 &0.0702 & 0.1103\\
& & ISB &0.0021 &0.0057 &0.0107 &0.0062 &0.0384 &0.0544 \\
& & IV &0.0193 &0.0153 &0.0147 &0.0464 &0.0318 &0.0559 \\
\hline
\end{tabular*}      \vspace*{3pt}
\end{table}

We also investigated how the additional terms in the modeling (\ref
{model}) affected the estimation precision when the true model was
given by
\[
g\bigl(m(\mathbf{x},\mathbf{z})\bigr)=x_1 f_1(z_1)
+ \cdots+ x_d f_d(z_d)
\]
for a set of covariates $(X_1, \ldots, X_d; Z_1, \ldots, Z_d)$. In the
latter model, each covariate appears in only one nonlinear interaction
term. For this, we estimated the following model:
%
\begin{eqnarray}
\label{simodel-xz} g\bigl(m(\bX,\bZ)\bigr)&=& f_{01}(Z_1) +
f_{02}(Z_2) + X_1\bigl(f_{11}(Z_1)+f_{12}(Z_2)
\bigr)
\nonumber
\\[-8pt]
\\[-8pt]
\nonumber
&&{} + X_2\bigl(f_{21}(Z_1)+f_{22}(Z_2)
\bigr),
\end{eqnarray}
where $f_{01}(z)=f_{02}(z)=0,   f_{11}(z)=\cos(2\pi z),
f_{12}(z)=0,  f_{21}(z)=0,  f_{22}(z)=\sin(2\pi z)$, and the link $g$
was the same as in the first example. The covariate~$X_1$ was a
discrete random variable having $\operatorname{Bernoulli}(0.5)$ distribution, $X_2$ was
the standard normal random variable and $Z_1$ and $Z_2$ were
$\operatorname{uniform}(0,1)$ random variables. The four covariates were independent.
The theoretically optimal bandwidths as defined at (\ref{optband}) were
$h_1^\mathrm{opt}=0.2405,   h_2^\mathrm{opt}=0.2469$ for $n=500$ and
$h_1^\mathrm{opt}=0.2093,   h_2^\mathrm{opt}=0.2149$ for $n=1000$, and we
used these in the simulation.

The main purpose of this additional simulation is to compare our
estimators based on the working model (\ref{simodel-xz}) with the
``oracle'' estimators which use the knowledge that
$f_{01}(z)=f_{02}(z)=f_{12}(z)=f_{21}(z)=0$. The system of updating
equations for the oracle estimators in our setting is given by (\ref
{backeq}) with the following modifications of $\hW_{jk}^{[s]}(z_j,z_k)$
and $\hF_j(\etam)(z_j)$: for $j \neq k$,
\begin{eqnarray*}
\hW_{jk}^{[s]}(z_j,z_k) &=& -
n^{-1}\sum_{i=1}^n
Q_2 \bigl(X_1^{i}{\hat{f}}_{11}^{[s]}(z_1)+X_2^i{
\hat{f}}_{22}^{[s]}(z_2), Y^i \bigr)
 X_j^i X_k^i
K_\bh\bigl(\bZ^i,\mathbf{z}\bigr),
\\
\hW_{jj}^{[s]}(z_j) &=& - \int n^{-1}
\sum_{i=1}^n Q_2
\bigl(X_1^{i}{\hat{f}}_{11}^{[s]}(z_1)+X_2^i{
\hat{f}}_{22}^{[s]}(z_2), Y^i \bigr)
\\
&&\hspace*{54pt}{}  \times \bigl(X_j^i\bigr)^2
K_\bh\bigl(\bZ^i,\mathbf{z}\bigr) \,d\mathbf{z}_{-j},
\\
\hF_j(\etam) (z_j)& =& \int n^{-1}\sum
_{i=1}^n Q_1
\bigl(X_1^{i}\eta_{11}(z_1)+X_2^i
\eta_{22}(z_2), Y^i \bigr)
 X_j^i K_\bh\bigl(
\bZ^i,\mathbf{z}\bigr) \,d\mathbf{z}_{-j},
\end{eqnarray*}
where $\bZ^i=(Z_1^i,Z_2^i)^\top$, $\mathbf{z}=(z_1,z_2)^\top$ and
$\etam(\mathbf{z}
)=(\eta_{11}(z_1), \eta_{22}(z_2))^\top$. Note that all these terms are
a scalar, not a matrix or a vector.
\begin{table}
\caption{Comparison of the smooth backfitting estimators
under the extended model (\protect\ref{simodel-xz})
and the corresponding oracle estimators}\label{tab2}
\begin{tabular*}{\textwidth}{@{\extracolsep{\fill}}lcccc@{}}
\hline
& & &$\bolds{f_{11}}$ &\multicolumn{1}{c@{}}{$\bolds{f_{22}}$}\\
\hline
Oracle &$n=500$ &IMSE &0.0680 &0.0639 \\
SBF& & ISB &0.0285 &0.0421 \\
& & IV &0.0395 &0.0218 \\[3pt]
&$n=1000$ &IMSE & 0.0400 & 0.0433 \\
& & ISB &0.0183 &0.0309 \\
& & IV &0.0216 &0.0124 \\[6pt]
SBF &$n=500$ &IMSE & 0.1057 & 0.0638 \\
based on (\ref{simodel-xz})& & ISB &0.0273 &0.0408 \\
& & IV &0.0784 &0.0230 \\[3pt]
&$n=1000$ &IMSE & 0.0627 & 0.0427 \\
& & ISB &0.0180 &0.0299 \\
& & IV &0.0447 &0.0128 \\
\hline
\end{tabular*}
\end{table}

Table~\ref{tab2} shows the results based on 500 datasets. For each of the
nonzero component functions, the table provides ISB, IV and IMSE. We
see that the smooth backfitting estimators perform fairly well in
comparison with their oracle versions. In particular, both have nearly
the same IMSE, ISB and IV for the estimation of
the second component function $f_{22}$. For estimating~$f_{11}$, the
smooth backfitting procedure with the extended model (\ref{simodel-xz})
gave almost the same bias as the oracle procedure, but a larger
variance than the latter. This may be expected since the former has the
additional component function $f_{12}$ in the estimation. This was not
the case with the estimation of $f_{22}$, however. The main reason is
that the variances of the estimators depend\vadjust{\goodbreak} highly on the design of the
regressor $X_2$. Recall that in parametric linear regression the
variance of the least squares estimator of a regression coefficient
gets smaller as the corresponding regressor is more variable. In our
setting, the variability of $X_2$ is four times as high as that of
$X_1$. This relatively high variability of $X_2$ alleviated the extra
sampling variability of the SBF estimator under the model
(\ref{simodel-xz}).

\begin{appendix}

\section*{Appendix: Technical details}\label{app}\vspace*{-1pt}
\setcounter{equation}{0}

\subsection{\texorpdfstring{Proof of Theorem \protect\ref{rateproposition}}{Proof of Theorem 1}}

The statement of Theorem~\ref{rateproposition} follows
immediately from the following lemma.\vspace*{-4pt}

\begin{lemma}\label{idenlem}
Under assumption \textup{(A0)}, it holds that there exist constants $0
< C_1 < C_2$ such that for two tuples $(\alpham,f_{jk}\dvtx  1 \le j
\le d, k \in I_j)$ and $(\alpham^*,f_{jk}^*\dvtx \break 1 \le j \le d, k \in
I_j)$ it holds that
%
\begin{eqnarray}
\label{claima1} && C_1 \int\bigl[g\bigl(m(\mathbf{x})\bigr)-g
\bigl(m^*(\mathbf{x})\bigr)\bigr]^2 P_{\bX}(d\mathbf{x})
\nonumber\\[-2pt]
&&\qquad\leq\bigl| \alpham-\alpham^*\bigr|_*^2 + \sum
_{j=1}^{d}\sum_{k \in
I_j}
\bigl[f_{jk}(x_{k}) - f^*_{jk}(x_{k})
\bigr]^2 p_{X_k}(x_k) \,d x_{k}
\\[-2pt]
\nonumber
&&\qquad \leq C_2 \int\bigl[g\bigl(m(\mathbf{x})\bigr)-g
\bigl(m^*(\mathbf{x})\bigr)\bigr]^2 P_{\bX
}(d\mathbf{x}).
\end{eqnarray}
Here, $p_{X_k}$ is the density of $X_k$, and
%
\begin{eqnarray}\quad
\label{struc-m} |\alpham|_*^2 &=& \int{ \Biggl( \sum
_{j=1}^d \alpha_j x_j +
\mathop{\sum_{ j
< k }}_{ j,k \in\mathcal{C}_0 } \alpha_{jk} x_j
x_k \Biggr)^2 P_{\bX
}(d\mathbf{x})},
\nonumber\\[-2pt]
\qquad\quad g\bigl(m(\mathbf{x})\bigr)& =& \sum_{j=1}^d
\alpha_j x_j + \mathop{\sum_{ j < k }}_{ j,k
\in
\mathcal{C}_0 }
\alpha_{jk} x_j x_k
+ {\tilde{\mathbf{x}}}_1^\top\bff_1(x_{r+1})
+ \cdots+ {\tilde{\mathbf{x}}}_p^\top\bff_p(x_{r+p}),
\\[-2pt]
g\bigl(m^*(\mathbf{x})\bigr) &=& \sum_{j=1}^d
\alpha^*_j x_j + \mathop{\sum_{ j < k
}}_{j,k \in\mathcal{C}_0 }
\alpha^*_{jk} x_j x_k
 + {\tilde{\mathbf{x}}}_1^\top\bff_1^*(x_{r+1})
+ \cdots+ {\tilde{\mathbf{x}}}_p^\top
\bff^*_p(x_{r+p}).\nonumber\vspace*{-3pt}
\end{eqnarray}
\end{lemma}

\begin{pf}
We only prove
the second inequality of (\ref{claima1}). The first one follows by
direct arguments. We first observe that because of assumption~(A0)
it holds that for constants $c_1, c_2> 0$,
\begin{eqnarray*}
&&\int\bigl[g\bigl(m(\mathbf{x})\bigr)-g\bigl(m^*(\mathbf{x})\bigr)
\bigr]^2 P_{\bX}(d\mathbf{x})
\\[-2pt]
&&\qquad \ge c_1 \int\bigl[g\bigl(m(\mathbf{x})\bigr)-g\bigl(m^*(
\mathbf{x})\bigr)\bigr]^2 \prod_{l=1}^D
P_{X_l}(dx_l)
\\[-2pt]
&&\qquad \ge c_2 \int\bigl[g\bigl(m(\mathbf{x})\bigr)-g\bigl(m^*(
\mathbf{x})\bigr)\bigr]^2 \prod_{l
\in
\cC}
w_l(x_l) \,dx_l \prod
_{l\notin\cC} P_{X_l}(dx_l).
\end{eqnarray*}
Denote by $\bI$ the right-hand side of the second inequality. Due
to the constraints of~(\ref{constraint1}) and
the fact that ${\tilde{\mathbf{x}}}_k$ does not include $x_{r+k}$, those terms
\begin{eqnarray*}
&\displaystyle\sum_{j=1}^d\bigl(\alpha_j-
\alpha^*_j\bigr) x_j + \mathop{\sum
_{ j < k}}_{ j,k
\in\cC_0 }\bigl(\alpha_{jk} - \alpha^*_{jk}
\bigr)x_j x_k,&
\\
&\displaystyle{\tilde{\mathbf{x}}}_1^\top\bigl(\bff_1(x_{r+1})-
\bff_1^*(x_{r+1})\bigr), \ldots, {\tilde{
\mathbf{x}}}_p^\top\bigl(\bff_p(x_{r+p})-
\bff_p^*(x_{r+p})\bigr)&
\end{eqnarray*}
are orthogonal in
$L_2(\mu)$, where $\mu$ is the product measure defined by
$\mu(d\mathbf{x})=\prod_{j\in\cC} w_j(x_j)   \,dx_j \prod_{j\notin
\cC}
P_{X_j}(dx_j)$. By this and by making use of (A0) again, we get
\begin{eqnarray*}
\bI&=& c_2 \int \Biggl[\sum_{j=1}^d
\bigl(\alpha_j-\alpha^*_j\bigr) x_j + \mathop{\sum
_{ j < k
}}_{
j,k \in\cC_0 }\bigl(\alpha_{jk} -
\alpha^*_{jk} \bigr)x_j x_k
\Biggr]^2 \prod_{l \in\cC} w_l(x_l)
\,dx_l \prod_{l \notin
\cC
}P_{X_l}(dx_l)
\\
&&{} + c_2 \sum_{k=1}^p \int
\bigl[{\tilde{\mathbf{x}}}_k^\top\bigl(
\bff_k(x_{r+k})- \bff^*_k(x_{r+k})
\bigr) \bigr]^2 \prod_{l\in\cC}
w_l (x_l) \,dx_l \prod
_{l \notin\cC
}P_{X_l}(dx_l)
\\
&\ge& c_3 \bigl|\alpham-\alpham^*\bigr|_*^2\\
&&{} + c_3\sum
_{k=1}^p \int \bigl[{\tilde{
\mathbf{x}}}_k^\top\bigl( \bff_k(x_{r+k})-
\bff^*_k(x_{r+k})\bigr) \bigr]^2
 P_{X_{r+k}}(dx_{r+k}) \prod
_{l \in I_{k}^*}^DP_{X_l}(dx_l)
\end{eqnarray*}
for some constants $c_3 > 0$ and where $I_{k}^*$ denotes the set
of indices of ${\tilde{\mathbf{x}}}_k$. The second inequality of (\ref{claima1})
now follows because the smallest eigenvalues of $ \int{\tilde{\mathbf{x}}}_k
{\tilde{\mathbf{x}}}_k^\top\prod_{l \in I_{k}^*}^D P_{X_l}(dx_l)$ can
be bounded
from below by a positive constant times the smallest eigenvalue of
$E[\tbX_k \tbX_k^\top] = \int{\tilde{\mathbf{x}}}_k {\tilde{\mathbf{x}}}_k^\top
P_{\tbX_k}(d{\tilde{\mathbf{x}}}_k)$, where $P_{\tbX_k}$ denotes the
distribution
of $\tbX_k$. These eigenvalues can be bounded away from zero by
assumption (A0).
\end{pf}

\subsection{\texorpdfstring{Proof of Corollaries \protect\ref{corratepen} and \protect\ref{corratesieve}}
{Proof of Corollaries 1 and 2}}

For the proof of these two corollaries, we apply Theorem
\ref{rateproposition}. We have to show that (\ref{estrate}) holds
with $\kappa_n= n^{l/(2l+1)}$ for the penalized least squares
estimator and the spline sieve estimator, respectively.

For the proof of Corollary~\ref{corratepen}, we apply Theorem 10.2
in \citet{geer}. As discussed in \citet{geer} the
statement of the theorem remains valid for errors with
subexponential tails if the entropy bounds hold for entropies
with bracketing. For the application of this theorem one needs
results on the entropy with bracketing for the class of functions
$m$ that fulfill (\ref{modelnormed}) with $ f_{jk}$\vadjust{\goodbreak} in Sobolev
classes. Because $g$ has an absolutely bounded derivative, Lemma
\ref{idenlem} implies that the well-known entropy conditions for
Sobolev classes carry over to the classes of functions $m$. This
proves Corollary~\ref{corratepen}.

For the proof of Corollary~\ref{corratesieve} we use Theorem 1 in
\citet{chen1998}. Compare also Theorem 10.11 in \citet{geer}.
Using the above entropy bound one can easily verify
Conditions A.1--A.4 in \citet{chen1998} with $l(\theta,
(\bX,Y))= (m(\bX)-Y)^2$, $\theta=(\alpham, f_{jk};1 \le j \le d, k
\in I_j)$ and $m$ as given at~(\ref{struc-m}). Note that
$l(\theta, (\bX,Y)) - l(\theta_0, (\bX,Y)) = (m(\bX)-m_0(\bX))^2-
2 (m(\bX)-m_0(\bX)) \varepsilon$, where $\theta_0=(\alpham_0,
f_{0jk};1 \le j \le d, k \in I_j)$ is the true tuple, $m_0$
denotes the true underlying regression function and
$\varepsilon=Y-m_0(\bX)$. To check the conditions compare also the
proof of Proposition 1 in the latter paper. In particular, their
condition A.4 holds with $s= 2l/(2l+1)$. This follows because for
two functions $g_1,g_2\dvtx [0,1] \to\mathbb{R}$ with $|D_z^lg_1 (z)|
\leq L$, $|D_z^lg_2 (z)| \leq L$ and $\int_0^1 (g_1(z) - g_2(z))^2
\,dz \leq\delta^2$, it holds that $|g_1 (z)-g_2 (z)| \leq2
(2L)^{1-c} \delta^{1-c}$ with $c= 2l(2l+1)^{-1}$; see Lemma 2 in
\citet{chen1998}. The necessary conditions are simplified
because we assume that the data are i.i.d.; see also Remark 1(b)
in \citet{chen1998}. To get $\varepsilon_n^{(2-s)/(\gamma-1)}
B_n \geq1$ at the end of their proof of Theorem~\ref{Thmdist} one needs that
$(2-s)/(\gamma-1) < s$ which is equivalent to $\gamma> 2+
l^{-1}$. One can check that their proof goes through with this
constraint. Thus it suffices for the i.i.d. case that $E |
\varepsilon|^{\gamma} < \infty$ holds for some $\gamma> 2+
l^{-1}$.

\subsection{Additional assumptions for kernel smoothing}

We assume the density of $\bX^{c}$ is supported on $[0,1]^p$.
Thus, the integration at (\ref{kernel}) is over $[0,1]$. We note
that, for the normalized kernel $K_{h_j}(u_j,z_j)$ introduced in
Section~\ref{theory}, $[2h_j, 1-2h_j]$ is the interior region for
$z_j$ that does not have a boundary effect. In addition to
assumption (A0) in Section~\ref{sec2}, we collect the conditions we use for
the theory in Sections~\ref{theory} and~\ref{lpsbf}.

\begin{longlist}[(A1)]
\item[(A1)] The quasi-likelihood function $Q(\mu,y)$ is three
times continuously differentiable with respect to $\mu$ for each
$y$ in its range, $Q_2(u, y) < 0$ for $u \in\R$ and~$y$ in its
range, the link function $g$ is three times continuously
differentiable, $V$ is twice continuously differentiable and the
conditional variance function $v(\mathbf{x})=\operatorname{ var}(Y|\bX=\mathbf{x})$ is
continuous in $\mathbf{x}^c=(x_{r+1}, \ldots, x_{r+p})^\top$ for each
$(x_1, \ldots, x_r)$. The densities $p_{X_j}$ for $r+1 \le j \le
r+p$ are bounded away from zero on $[0,1]$. The function $V$ and
the derivative $g^{\prime}$ are bounded away from zero. The
higher-order derivatives $g^{\prime\prime}$ and $g^{\prime\prime
\prime}$ are bounded. The weight function $w$ is continuously
differentiable and fulfills $w(0)=w(1)=0$.
\item[(A2)] The partial
derivatives $\partial p_{\bX}(\mathbf{x})/\partial\mathbf{x}^c$ of the joint
density function $p_{\bX}$ exist and are continuous in $\mathbf{x}^c$ for
all $(x_1, \ldots, x_r)$.
\item[(A3)] The components of $\bff_j$
are twice continuously differentiable.
\item[(A4)] $E|Y|^\alpha<
\infty$ for some $\alpha> 5/2$.
\item[(A5)] The kernel function
$K$ is bounded, symmetric about zero, has compact support, say
$[-1,1]$, and is Lipschitz continuous. The bandwidths $h_j$ depend
on the sample size $n$ and satisfy $n^{1/5}h_j \rightarrow c_j$ as
$n \rightarrow\infty$ for some constants $0<c_j<\infty$.
\end{longlist}

\subsection{\texorpdfstring{Preliminaries for the proofs of theorems \protect\ref{Thmrate}--\protect\ref{Thmalgo}}
{Preliminaries for the proofs of theorems 2--4}}\label{prelim}

The population versions of $\hbF_j$ are defined by
\[
\bF_j(\etam) (z_j)= \int E \bigl[Q_1
\bigl(\tbX^\top\etam(\mathbf{z}), m(\bX) \bigr) \tbX_j |
\bX^c=\mathbf{z} \bigr] p_{\bX^c} (\mathbf{z}) \,d
\mathbf{z}_{-j}.
\]
For the empirical versions of $\bW_{jk}(\mathbf{z};\etam)$,
$\bW_j(\mathbf{z};\etam)$ and $\bW(\mathbf{z};\etam)$ introduced in Section~\ref{theory},
we define
\[
\hbW_{jk}(\mathbf{z};\etam) = - n^{-1}\sum
_{i=1}^n Q_2 \bigl(\tbX^{i\top}
\etam(\mathbf{z}), Y^i \bigr)\tbX_j^i
\tbX_k^{i\top} K_\bh\bigl(\bX^{c,i},
\mathbf{z}\bigr),
\]
and then define $\hbW_{j}(\mathbf{z};\etam)$ and $\hbW(\mathbf{z};\etam)$
in the
same way as we define $\bW_j(\mathbf{z};\etam)$ and $\bW(\mathbf{z};\etam)$,
respectively. We write $\hbW_{jk}(\mathbf{z})=\hbW_{jk}(\mathbf{z};\bff)$ in
case the true $\bff$ enters into the place of $\etam$.

For a tuple $\deltam\in\cH(\bW)$, let $\hbF_j'(\etam)(\deltam)$
denote the Fr\'echet differential of $\hbF_j$ at $\etam$ to the
direction of $\deltam$. Then
\[
\hbF_j'(\etam) (\deltam) (z_j) = -\int
\hbW_j(\mathbf{z};\etam) \deltam(\mathbf{z}) \,d\mathbf{z}_{-j}.
\]
The second term on the right-hand side of (\ref{NR}) is simply
$\hbF_j'({\hat{\bff}}^{[s-1]})({\hat{\bff}}- {\hat{\bff
}}^{[s-1]})(z_j)$. The population
versions of $\hbF_j'(\etam)$ are defined by
$\bF_j'(\etam)(\deltam)(z_j) = -\int\bW_j(\mathbf{z};\etam)
\deltam(\mathbf{z})  \,d\mathbf{z}_{-j}$. Define a linear operator $\hbF'(\etam)$
by
\[
\hbF'(\etam) (\deltam) = \bigl(\bigl(\hbF_1'(
\etam) (\deltam)\bigr)^\top, \ldots, \bigl(\hbF_p'(
\etam) (\deltam)\bigr)^\top \bigr)^\top.
\]
Likewise, define $\bF'(\etam)$ from $\bF_j'(\etam)$. In the proofs
below, we use $\bff=(\bff_1^\top,
\ldots, \bff_p^\top)^\top$ to denote the true vector of univariate
functions.

\subsection{\texorpdfstring{Proof of Theorem~\protect\ref{Thmrate}}{Proof of Theorem 2}}

In addition to $\|\cdot\|_\#$ introduced in Section~\ref{theory},
we consider two other norms. Let $\|\cdot\|_2$ be the
$L_2(p_{\bX^c})$-norm defined by $\|\etam\|_2^2 = \int
|\etam(\mathbf{z})|^2   p_{\bX^c}(\mathbf{z})  \,d\mathbf{z}$. Define
$\|\etam\|_{\infty} = \max\{\sup_{2h_1 \le z_1 \le1-2h_1}|\etam_1(z_1)|, \ldots, \break
\sup_{2h_p \le z_p \le1-2h_p}|\etam_p(z_p)|\}$,
where $|\cdot|$ denotes the Euclidean norm. As in Section~\ref{theory},
we write $\bW(\mathbf{z})=\bW(\mathbf{z};\bff)$, $\bW_{jj}(z_j)=\bW_{jj}(z_j;\bff)$,
etc., for the true tuple $\bff$. For a linear operator $\cF$ that maps
$\cH(\bW)$ to $\cH(\bW)$, let $\|\cF\|_\mathrm{ op}$ denote its
operator-norm defined by $\|\cF\|_\mathrm{ op} = \sup_{\|{\footnotesize
\deltam}\| =1} \|\cF(\deltam)\|$. Here and below, if not specified,
$\|
\cdot\|$ is either $\|\cdot\|_2$ or $\|\cdot\|_\infty$. We prove
%
\begin{eqnarray}
\label{proof1} P \bigl(\hbF'(\bff) \mbox{ is invertible and }\bigl\|
\hbF'(\bff )^{-1}\bigr\|_\mathrm{ op} \le C_1
\bigr)& \rightarrow&1,
\\
P \bigl(\bigl\|\hbF'(\etam)-\hbF'\bigl(
\etam'\bigr)\bigr\|_\mathrm{ op} \le C_2 \bigl\|\etam-
\etam'\bigr\| \mbox{ for all } \etam, \etam' \in
B_r(\bff) \bigr) &\rightarrow& 1, \label{proof2}
\end{eqnarray}
for some constants $r,C_1,C_2>0$, where $B_r(\bff)$ is a ball centered
on $\bff$ with radius $r$.
Then, the theorem follows from Newton--Kantorovich theorem\vadjust{\goodbreak} [see,~e.g., \citet{deimling}] since $\|\hbF(\bff)\|_2 = O_p(n^{-2/5})$ and
$\|
\hbF(\bff)\|_\infty= \break O_p(n^{-2/5}\sqrt{\log n})$.

By the standard techniques of kernel smoothing, one can show that,
uniformly for $\mathbf{z}\in[0,1]^p$,
$\hbW_{jk}(z_j, z_k)\equiv\int\hbW_{jk}(\mathbf{z}) \,d\mathbf{z}_{-(j,k)}$
converges to $ \bW_{jk}(z_j, z_k)$ and $ \hbW_{jj}(z_j) \equiv
\int\hbW_{jj}(\mathbf{z}) \,d\mathbf{z}_{-j}$ to $\bW_{jj}(z_j)$, for $1
\le j
\neq k \le p$. This gives the
uniform convergence of $\hbF'(\bff)(\deltam)$ to
$\bF'(\bff)(\deltam)$ over $\deltam$ such that $\|\deltam\| \le
R$, where $R>0$ is an arbitrary positive real number. Thus, to
prove (\ref{proof1}) it suffices to show that $\bF'(\bff)$ is
invertible and has a bounded inverse. For this claim we first show
that the map $\bF'(\bff)\dvtx  \cH(\bW) \rightarrow\cH(\bW)$ is
one-to-one. Suppose that $\bF'(\bff)(\deltam) =\bzero$ for some
$\deltam\in\cH(\bW)$. We have to show that $\deltam=\bzero$.
From $\bF'(\bff)(\deltam) =\bzero$ we get that
$\deltam_{j}(z_j)^\top\int\bW_j(\mathbf{z}) \deltam(\mathbf{z})  \,d\mathbf{z}_{-j} =
0$ for all $1 \le j \le p$. This implies
\begin{eqnarray*}
0&=& \sum_{j=1}^p \int \biggl[\int
\deltam_{j}(z_j)^\top\int\bW_j(
\mathbf{z}) \deltam(\mathbf{z}) \,d\mathbf{z}_{-j} \biggr]
\,dz_j
\\
&\ge& \sum_{j=1}^p \sum
_{k=1}^p \int\bigl({\tilde{\mathbf{x}}}_j^\top
\deltam_{j}(x_{r+j})\bigr) \bigl({\tilde{
\mathbf{x}}}_k^\top \deltam_{k}(x_{r+k})
\bigr)P_{\bX}(d\mathbf{x})
\\
&=& c \int \Biggl(\sum_{j=1}^p {\tilde{
\mathbf{x}}}_j^\top\deltam_{j}(x_{r+j})
\Biggr)^2 P_{\bX
}(d\mathbf{x})
\end{eqnarray*}
for some positive constant $c > 0$. Here, for the inequality we
used that $V(u)g'(u)^2$ is bounded from above for $u$ in any
compact set. Applying the arguments in the proof of Lemma~\ref{idenlem}, we
get that $\deltam= \bzero $ a.s. This proves the claim that the
map $\bF'(\bff)\dvtx  \cH(\bW) \rightarrow\cH(\bW)$ is one-to-one.

Next, using the fact that $\langle\bF'(\bff)(\deltam), \etam
\rangle_\# =\langle\deltam, \bF'(\bff)(\etam) \rangle_\#$ for all
$\deltam,\etam\in\cH(\bW)$, one can show that $\bF'(\bff)$ is
onto. Thus, $\bF'(\bff)$ is invertible. To verify that
$\bF'(\bff)^{-1}$ is bounded, it suffices to prove that the
bijective linear operator $\bF'(\bff)$ is bounded, due to the bounded
inverse theorem. From the assumption~(A1) and the fact that
support of the density of $\bX$ is bounded, we get
\begin{eqnarray*}
\bigl\|\bF'(\bff) (\deltam)\bigr\|_\#^2 &=& \int\bigl|
\bF'(\bff) (\deltam) (\mathbf{z})\bigr|^2 p_{\bX^c}(
\mathbf{z}) \,d\mathbf{z} \le C_3 \|\deltam\|_2^2,
\\
\bigl\|\bF'(\bff) (\deltam)\bigr\|_\infty&\le& C_4 \|
\deltam\|_\infty
\end{eqnarray*}
for some constants $C_3, C_4>0$. This concludes that $\bF'(\bff)$
is bounded in both of the norms $\|\cdot\|_2$ and
$\|\cdot\|_\infty$.

The claim (\ref{proof2}) holds since, for any given $r>0$, there
exists a constant $C_5>0$ such that, with probability tending to
one, $\|\hbF'(\etam)(\deltam)-\hbF'(\etam')(\deltam)\| \le C_5
\|\etam- \etam'\| \cdot\|\deltam\|$ for all $\etam,\etam' \in
B_r(\bff)$.

\subsection{\texorpdfstring{Proof of Theorem~\protect\ref{Thmdist}}{Proof of Theorem 3}}

Let $\cdelm$ denote a solution of the following equations:
%
\begin{eqnarray}
\label{backeq1} \qquad\cdelm_{j}(z_j)& =&
\tdelm_j(z_j)
\nonumber
\\[-8pt]
\\[-8pt]
\nonumber
&&{}- \sum_{k\neq j}
\int \hbW_{jj}(z_j)^{-1}\hbW_{jk}(z_j,z_k)
\cdelm_{k}(z_k) \,dz_k,\qquad
 1 \le j \le p,
\end{eqnarray}
where $\tdelm_j(z_j) = \hbW_{jj}(z_j)^{-1} \hbF_j(\bff)(z_j)$. We
first remark that $\cdelm$ exists and is unique, with probability
tending to one. Define $\tilde{\bff}=\bff+ \cdelm$. We claim
%
\begin{equation}
\label{proof15} \|\tilde{\bff}-\bff\|_\infty= O_p
\bigl(n^{-2/5}\sqrt{\log n}\bigr).
\end{equation}
Define $\tilde\bF_j$ for $1 \le j \le p$ by
\begin{eqnarray*}
\tilde\bF_j (\etam) (z_j) &=& \hat\bF_j (
\bff) (z_j)\\
&&{}+ \sum_{k=1}^p
\int n^{-1}\sum_{i=1}^n
Q_2 \bigl(\tbX^{i\top}\bff(\mathbf{z}), Y^i \bigr)
\\
&&\hspace*{74pt} {}\times \tbX_j^i \tbX_k^{i\top}
\bigl[\etam_{k}(z_k)-\bff_{k}(z_k)
\bigr]K_\bh\bigl(\bX^{c,i},\mathbf{z}\bigr) \,d
\mathbf{z}_{-j}.
\end{eqnarray*}
Note that $\tilde\bff$ is the solution of the system of
equations $\tilde\bF_j (\etam)=\bzero_j,   1 \le j \le p$ by the
definitions of $\cdelm_j$ and $\tdelm_j$. Thus, the claim
(\ref{proof15}) ensures that $\hbF_j(\tilde\bff)(z_j) =
\tilde\bF_j (\tilde\bff)(z_j) + o_p(n^{-2/5}) = o_p(n^{-2/5})$,
uniformly for $z_j \in[2h_j, 1-2h_j]$. Also, (\ref{proof15}) and
(\ref{proof2}) give (\ref{proof1}) with $\bff$ being replaced by
$\tilde\bff$. This establishes $\|{\hat{\bff}}-\tilde{\bff}\|_\infty=
o_p(n^{-2/5})$. The theorem follows since
\[
n^{2/5}\cdelm(\mathbf{z}) \mathop{\rightarrow}^{d} N \bigl(
\bigl(\betam_1(z_1)^\top, \ldots,
\betam_p(z_p)^\top\bigr)^\top, \operatorname{
diag}\bigl(\Sigmam_j(z_j)\bigr) \bigr),
\]
the latter being proved similarly as in the proof of Theorem~2 of
Lee, Mammen and Park (\citeyear{lee2012}).

It remains to prove (\ref{proof15}). The Fr\'echet differential of
$\tilde\bF_j$ at
$\etam$ to the direction $\deltam$, which we denote by
$\tilde\bF_j'(\etam)(\deltam)$, does not depend on $\etam$ since
$\tilde\bF_j(\etam)$ is linear in $\etam$. In fact
$\tilde\bF_j'(\etam)(\deltam) = \hbF_j'(\bff)(\deltam)$ for all
$\etam$. This means that $\tilde\bF(\bff)= \hbF(\bff)$ and
$\tilde\bF'(\etam)= \hbF'(\bff)$ for all $\etam$, so that
(\ref{proof1}) and (\ref{proof2}) are valid for $\tilde\bF$. As in the
proof of Theorem~\ref{Thmrate},
this implies (\ref{proof15}).

\subsection{\texorpdfstring{Proof of Theorem~\protect\ref{Thmalgo}}{Proof of Theorem 4}}
An application of Newton--Kantorovich theorem gives the first part
of the theorem. For the proof of the second part of the theorem, we
rewrite a full
cycle of the iteration step in (\ref{innerit}) as
$\hdelm^{[s,\ell]} = \hat{\deltam}_+^{[s-1]} + \hat A^{[s-1]}
\hdelm^{[s,\ell-1]}$
with appropriate definitions of $\hat{\deltam}_+^{[s-1]}$ and $\hat
A^{[s-1]}$.
Note that $\hat{\deltam}_+^{[s-1]}$ differs from the tuple with
elements $\tdelm_j^{[s-1]}$.
Also, we can write a~full cycle of the
iteration step for solving (\ref{backeq1}) as $\check{\deltam
}^{[\ell]}
= \check{\deltam}_+ + \check A\check{\deltam}^{[\ell-1]}$ with
appropriate\vadjust{\goodbreak}
definitions of $\check{\deltam}_+$ and $ \check A$. Finally, we can
write $\deltam^{[\ell]} = \deltam_+ + A\deltam^{[\ell-1]}$ with
appropriate definitions of $\deltam_+$ and~$A$ for a full cycle of
\begin{eqnarray}
\deltam_{j}^{[\ell]}(z_j) =
\bW_{jj}(z_j)^{-1} \bF_j(\bff)
(z_j) - \sum_{k\neq j} \int\bigl[
\bW_{jj}(z_j)\bigr]^{-1} \bW_{jk}(z_j,z_k)
\deltam_{k}^{[\ell-1]}(z_k)
\,dz_k,\nonumber\\
\eqntext{ 1 \le j \le p.}
\end{eqnarray}
For the convergence of the last iteration, we note that the
projection operators ${\pi}_{kj}\dvtx  \cH_k (\bW) \rightarrow
\cH_j(\bW)$ for all $1 \le j \neq k \le p$ are Hilbert--Schmidt,
where $\cH_k(\bW)$ is a subspace of $\cH(\bW)$ such that $\etam
\in\cH_k(\bW)$ if and only if $\etam_l =\bzero$ for all $l \neq
k$, and elements of $\etam_k$ with the configuration at
(\ref{repre}) satisfy the constraints (\ref{constraint3}). This implies
there exist constants $C_0$ and
$0<\rho_0<1$, with
%
\begin{equation}
\label{proof3-5} \int\bigl[\deltam^{[\ell]}(\mathbf{z}) -
\deltam^{[\infty]}(\mathbf{z})\bigr]^\top \bW(\mathbf{z})\bigl[
\deltam^{[\ell]}(\mathbf{z}) - \deltam^{[\infty]}(\mathbf{z})\bigr] \,d
\mathbf{z} \le {C_0}\rho_0^{2\ell}
\end{equation}
for some limiting function $\deltam^{[\infty]}$.

We now apply $\|\check A - A\| = o_p(1)$, $\|\check{\deltam}_+ -
\deltam_+ \| =o_p(1)$ and $\sup_s\|\hat A^{[s]} -\check A \| \leq
c$, $\sup_s \|\hat{\deltam}_+^{[s]} - \check{\deltam}_+\| \leq c$
with probability tending to one, for some constant $c>0$ that can
be made as small as we like by choosing $\tau$ small enough. This and
equation~(\ref{proof3-5}) implies that for
some constants $C_*$ and $0<\rho_*<1$,
\[
\int\bigl[\hdelm^{[s,\ell]}(\mathbf{z}) - \hdelm^{[s,\infty]}(\mathbf{z})
\bigr]^\top \bW(\mathbf{z})\bigl[\hdelm^{[s,\ell]}(\mathbf{z}) -
\hdelm^{[s,\infty]}(\mathbf{z})\bigr] \,d\mathbf{z} \le {C_*}\rho_*^{2\ell},
\]
with probability tending to one. This completes the proof of
Theorem~\ref{Thmalgo}.
\end{appendix}

%


\printaddresses

\end{document}